\documentclass[11pt,a4paper]{article}
\usepackage[left=2.2cm,top=3.8cm,bottom=3cm,right=2.2cm]{geometry}

\usepackage[utf8x]{inputenc}
\usepackage{amsmath,amssymb}
\usepackage{xcolor}
\usepackage{mathtools}
\usepackage{tabulary}
\usepackage{amsthm}
\usepackage{latexsym}
\usepackage[mathcal]{eucal}
\usepackage{amscd}
\usepackage{graphicx}
\usepackage{amsfonts}
\usepackage{amscd}
\usepackage{MnSymbol}
\usepackage{setspace}
\newtheorem{definition}{Definition}[section]
\usepackage[english]{babel}

\newtheorem{remark}{Remark}[section]
\newtheorem{theorem}{Theorem}[section]
\newtheorem{lemma}{Lemma}[section]
\newtheorem{proposition}{Proposition}[section]
\newtheorem{corollary}{Corollary}[section]
\newtheorem{example}{Example}[section]

\theoremstyle{remark}

\usepackage{graphicx}
\date{}

\newcommand{\bt}{\begin{theorem}}
\newcommand{\et}{\end{theorem}}
\newcommand{\bl}{\begin{lemma}}
\newcommand{\el}{\end{lemma}}
\newcommand{\bexc}{\begin{exercise}}
\newcommand{\eexc}{\end{exercise}}
\newcommand{\bpr}{\begin{proposition}}
\newcommand{\epr}{\end{proposition}}
\newcommand{\bex}{\begin{example}}
\newcommand{\eex}{\end{example}}
\newcommand{\bc}{\begin{corollary}}
\newcommand{\ec}{\end{corollary}}
\newcommand{\bo}{\begin{proof}}
\newcommand{\eo}{\end{proof}}
\newcommand{\bd}{\begin{definition}}
\newcommand{\ed}{\end{definition}}
\newcommand{\br}{\begin{remark}}
\newcommand{\er}{\end{remark}}
\newcommand{\be}{\begin{enumerate}}
\newcommand{\ee}{\end{enumerate}}

\newcommand{\K}{\mathcal{K}}
\newcommand{\M}{\mathcal{M}}

\newcommand{\fS}{\mathfrak{S}}
\newcommand{\sft}{\mathfrak{sft}}

\def \ep {\epsilon}

\def \s {\sigma}
\def \o {\omega}

\newcommand{\A}{\mathcal{A}}
\newcommand{\cL}{\mathcal{L}}
\newcommand{\cR}{\mathcal{R}}

\newcommand{\cF}{\mathcal{F}}

\newcommand{\cO}{\mathcal{O}}

\newcommand{\Z}{\mathbb{Z}}

\newcommand{\N}{\mathbb{N}}

\begin{document}
\title{Finiteness in Polygonal Billiards on Hyperbolic Plane}
\author{Anima Nagar and Pradeep Singh\\
Department of Mathematics, Indian Institute of Technology Delhi,\\
Hauz Khas, New Delhi 110016, INDIA}


\maketitle

\vspace{1cm}

 \maketitle

\begin{abstract}
\textsc{J. Hadamard} studied the geometric properties of geodesic flows on surfaces of negative curvature, thus initiating ``Symbolic Dynamics". 

In this article, we follow the same geometric approach to study the geodesic trajectories of billiards in ``rational polygons" on the hyperbolic plane. We particularly show that the billiard dynamics  resulting thus are just `Subshifts of Finite Type' or their dense subsets. We further show that `Subshifts of Finite Type' play a central role in subshift dynamics and while discussing the topological structure of the space of all subshifts, we demonstrate  that they approximate any shift dynamics. 

\end{abstract}

\vspace{6cm}

\thanks{\emph{keywords:} hyperbolic plane, polygonal billiards, pointed geodesics, subshifts of finite type, Hausdorff metric, space of all subshifts.   }


\newpage

\tableofcontents

\section{Introduction}

`Mathematical billiards' describe the motion of a  point mass in a domain with elastic reflections from the boundary, and occur naturally in  many problems in science.  The billiards problem has typically been studied in planar domains.

A \emph{billiard} in a `domain' $\Pi$ in the Euclidean plane is defined as a dynamical system described by the motion of a point-particle within $\Pi$ along the straight lines with specular reflections from the boundary $\partial \Pi$.    A domain is generally taken to be a subset of the plane that is compact with a piecewise smooth boundary. We refer to the point-particle under consideration as the \emph{billiard ball}, the path followed within $\Pi$ as the \emph{billiard trajectory} and the respective domain is called the \emph{billiard table}.         This simple to describe mathematical system captures the essence of  `chaotic dynamics'. The resultant `dynamics' is heavily dependent upon the nature of the boundary of the billiard table and the position as well as the orientation of the onset of the billiard. The curvature of the part of the boundary being hit decides whether two parallelly  launched billiard balls will come out parallel, grow apart or will cross their paths. The overall trajectory is also determined by the relative placement of the pieces of the boundary with respect to each other. In this setting, the problem of billiards has been studied classically in \cite{bunimovich,gutkin3,gutkin,gutkin2,tabachnikov}. Although this problem in itself is still ripe with many interesting open problems \cite{gutkin4,schwartz1,schwartz2}. 

\bigskip

\emph{Symbolic dynamics} is the dynamical study
of the shift automorphism on the space of bi-infinite sequences on some set of symbols. The beginning of symbolic dynamics can be traced back to \textsc{J. Hadamard}, when in 1898 he studied geodesics on surfaces of negative curvature   from a geometric point of view. (See J. Hadamard, {\it Les surfaces \`a courbures oppos\'ees et leurs lignes g\'eod\'esiques}, J.Math.Pures Appl. 4 (1898), 27-73.
) Though a systematic study of symbolic dynamics is said to have started with the work of \textsc{Marston Morse and Gustav Hedlund} \cite{morse}. 

\bigskip

Hadamard can also be credited to first study \emph{Hadamard's billiards} which also counts as the first example proved to be possessing deterministic chaos. Along the same period of time, \textsc{G. D. Birkhoff} working under the ambit of relativity and discrete dynamical systems was quite interested in billiards, while studying the \emph{three-body problem}. Many initial results on such geodesic flows were compiled by \textsc{ Gustov Hedlund} \cite{hedlund}. Further work on Hadamard's billiards was carried forward by \textsc{M. Gutzwiller} in 1980s, see \cite{gutzwiller2,gutzwiller3, gutzwiller}. The study of polygonal billiards in hyperbolic plane besides its theoretical interest, appears in `General Relativity' in an extended form of polyhedral billiards in a \emph{hyperbolic space} (an n-dimensional manifold with constant negative `Riemannian curvature'), see \cite{castle, gutzwiller}. Another approach can be seen in \cite{skatok3}.

\bigskip

For planar billiards in the \emph{Euclidean plane}, it has been observed that  the route of dealing with   `Symbolic Dynamics' is inherently ambiguous as this setting lacks in providing a one-to-one correspondence between the billiard trajectories and the corresponding natural codes generated by collecting the labels of the sides of polygons being hit in order.

\bigskip

This problem can be handled better in the \emph{hyperbolic plane}.

\bigskip

A hyperbolic plane is defined as a 2-dimensional manifold with a constant negative Gaussian curvature. Since the hyperbolic plane cannot be embedded in a 3-dimensional Euclidean space, we are forced to work with various models of the hyperbolic plane, most common being the \emph{Poincar\'e half plane model $\mathbb{H}$} and the \emph{Poincar\'e disc model $\mathbb{D}$}.	We prefer the disc model for our work as it is Euclidean compact and this lets us follow the trajectories to each boundary which is the essence of traditional billiards. This is not the case with the half plane model.

\bigskip

   In this setting, we carry out our investigations with a class of polygons(being a subset of $\mathbb{D}$) whose boundary comprises of finitely many geodesic segments, is piecewise smooth and intersects in vertices at angles either $0$ or those that divide $\pi$ into integer parts.  This `rationality' at vertices lying in $\mathbb{D}$ is necessitated by the general technique that we follow which intrinsically takes the \emph{tiling} of billiard tables (i.e. tessalations of the ambient hyperbolic plane) into account. For this class of polygons as billiard tables, we show that the associated billiard trajectories correspond to countably many \emph{bi-infinite sequences}. The work on these lines was initiated in \cite{ullmo,ullmo2} and later pursued in \cite{castle}. 
   
   \bigskip

We achieve uniqueness in the `coding' by taking a billiard trajectory and breaking it into several \emph{`pointed geodesics'}, each described by a \emph{`base arc'}. These `pointed geodesics' can be regarded as compact subsets of $\mathbb{D}$ and form a space endowed with the Hausdorff topology. Under this splitting procedure, we establish a one-to-one correspondence between the space of `pointed geodesics' and the corresponding space of  bi-infinite sequences which we call as the associated \emph{shift space}.  The `dynamics' associated with these billiards is independent of the position and orientation of the polygon in $\mathbb{D}$, that allows us to choose any one of the isometric images appearing in the tiling as the `fundamental polygon'.

We introduce a metric $d_{\mathbb{G}}$ on the space ${\mathbb{G}}$ of pointed geodesics which turns out to be equivalent to the classical `Hausdorff metric'. Under this metric and the transformation $\tau$ provided by the `bounce map'(the map that describes the specular reflection of the billiard ball at each hit with the boundary), we consider $( \mathbb{G},\tau)$ to be a topological dynamical system. Further, collecting the symbols on the sides being hit by the billiard ball, we establish a `coding' for the trajectory. Thus, we get a `shift invariant' set $X$ of codes. We prove that $ {X}$ is a `shift of finite type (SFT)' that is not a `full k-shift' or its dense subset. We establish a conjugacy between $( \mathbb{G},\tau)$ and $( X, \sigma )$, thereby instituting an explicit route for studying the geometrical properties of billiards on this class of polygons via the symbolic dynamics on the corresponding space of  bi-infinite sequences.

\bigskip

At this point we would like to mention that such  representations of geodesic flows by symbolic systems have also been studied by \textsc{Roy Adler} and \textsc{Leopold Flatto} \cite{adler-flatto} and \textsc{Caroline Series} \cite{series1, series2}. One could also look into \cite{skatok2} for a more recent exposition. 

The two dimensional geometry that we consider is similar to  \cite{adler-flatto}, that also looks into the relation between geodesic flows on a compact surface of constant negative curvature and their associated symbolic dynamics. Though, their approach is to look into geodesic flows on a compact surface $S = \mathbb{D}/\Gamma$ of genus $g \geq 2$ determined by the acting Fuchsian group $\Gamma$, as  for such a surface a bounded fundamental domain $F$ is guaranteed to exist with $8g-4$ sides.  The \emph{cutting sequences} that they refer to are the\emph{ pointed geodesics} in our case, though our concept of \emph{polygon} is very distinct from the \emph{fundamental domain} that they consider. Moreover, in \cite{adler-flatto, series1, series2} the discussions follow on the measure-theoretic lines whereas our constructions and proofs run on topological arguments.

\smallskip

It is noted here that our current problem has a setting that is different from the setting in \cite{adler-flatto, series1, series2}. The differences arise on several counts: the foremost being that the ambient space  considered in \cite{adler-flatto, series1, series2} is a compact surface of genus $g \geq 2$ whereas in our case, we consider the whole $\mathbb{D}$ in 	the background, thus allowing the billiard trajectories to lie on a non-compact surface. In \cite{adler-flatto, series1, series2} the authors start by choosing a \emph{Fuchsian group} $\Gamma$ such that $\mathbb{D}/ \Gamma $ is a compact surface of genus $g \geq 2$ and then for such a surface a bounded fundamental domain $F$ is guaranteed to exist with $8g-4$ sides. On the other hand, we work with a semi-ideal or ideal or compact rational polygons $\Pi$ on $\mathbb{D}$ which can even have  an odd number of sides and assume no identification of sides of $\Pi$.

In \cite{adler-flatto}, the authors have studied connection between the Geodesic flows and interval maps. For $(z,v) \in T^1 \mathbb{D}$, let $(\gamma(t))_{t \in \mathbb{R}}$ be the arclength parameterization of the unique geodesic passing through $z \in \mathbb{D}$ with the tangent line at $z$ in the direction $v$ and satisfying $\gamma(0) = z$ and $\gamma'(0) = v$. The \emph{Geodesic flow} on $T^1 \mathbb{D} $ is the map $\phi$ from $\mathbb{R} \times T^1 \mathbb{D}$ into $T^1 \mathbb{D}$ defined by $\phi(t, (z,v)) = (\gamma(t), \gamma'(t))$.

\begin{figure}[h!]
	\centering
	\includegraphics[width=10cm,height=9cm]{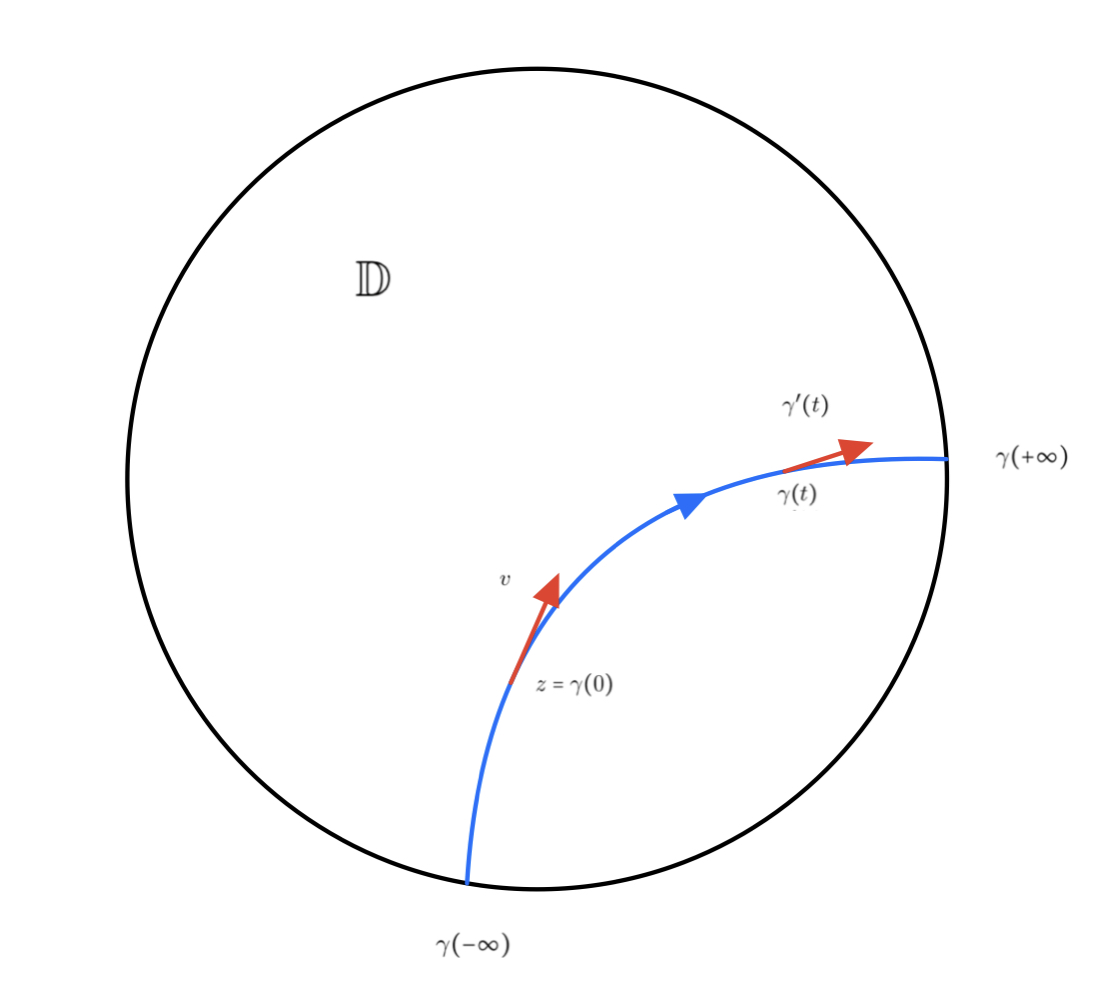}
	\caption{Geodesic flow on $\mathbb{D}$}
\end{figure}

The geodesic flow on the quotient of $T^1 \mathbb{D}$ by a Fuchsian group $\Gamma$ is defined in the natural way via the projection map $\pi : T^1 \mathbb{D} \rightarrow T^1(\mathbb{D}/ \Gamma)$. Thus, to link the continuous time geodesic flow with the discrete interval maps, various reductions are required. This is done by choosing a suitable \emph{cross section} and the corresponding \emph{cross section map}. A cross section is roughly considered as a subset $C$ of $T^1(\mathbb{D}/ \Gamma)$ which intersects the geodesic flow repeatedly in past and future. There are usually multiple choices available for the cross section. 

In our problem we choose the cross section explicitly to be $T^1(\partial \Pi)$ or under the projection map, simply $\partial \Pi$. Intuitively, the part of the billiard trajectory between two consecutive hits is rendered irrelevant. We lift the billiard trajectories to the corresponding pointed geodesics on route to establishing conjucacy between the corresponding space of pointed geodesics and the associated space of codes. 

On the contrary in \cite{adler-flatto}, the authors choose a cross section $C$ for which a second reduction is possible, which is given by a one-dimensional factor map. Through these two reductions, the relationship between the geodesic flows and interval maps is further discussed. They then relate these interval maps with symbolic systems via the \emph{Markov partitions}. This line of separation between the two is crucial as we are able to discuss the case of billiards even when a vertex of polygon $\Pi$ goes to $\partial \mathbb{D}$ whereas in \cite{adler-flatto} the inherent restriction of $\mathbb{D} / \Gamma$ being compact is maintained throughout. The existence of the fundamental region $F$ is essentially based on this. Our model is fairly simple.

\medskip

 The underlying constructions are different in nature and this results in a difference in the nature of the associated symbolic sequences that they obtain in \cite{adler-flatto} from what we obtain here.

\bigskip

Our results here motivate us to look into the larger picture of shift dynamics. We study the space of all closed shift invariant sets of  bi-infinite sequences on finite symbols, given the Hausdorff topology. We particularly demonstrate that `shifts of finite type (SFTs)' form a dense subset of such a space. This gives us convergence in subshifts and we further study the dynamical properties induced by such a convergence.

\bigskip

In Section 2. we discuss some basic theory thus introducing our definitions and notations. We study our billiard dynamics and derive the associated symbolic dynamics and their properties in Section 3. In Section 4, we study convergence of the dynamics in  subshifts and polygonal billiards.

\section{Preliminaries}
In this section, we lay down some basic notions for  later usage, thus establishing the notation we will henceforth use.

\subsection{Some Topological Dynamics}

  A discrete dynamical system $(X,f)$ consists of a continuous self-map $f$ on a  metric space $(X,d)$.

  \bigskip

  The \emph{orbit} of a point
$x \in X$ is the set $\mathcal{O}(x) = \{ x, f(x),f^2(x), \dots \}$.  Here $f^n$ stands for the
$n-$fold self-composition of $f$. We note that orbits are invariant sets i.e. $f(\cO(x)) \subseteq \cO(x)$ for all $x \in X$. The basic study in dynamics is to study the asymptotic behaviour of orbits of all $x \in X$.

The point $x$ is \emph{periodic} if there exists $n \in \N$ such that $f^n(x) = x$.  The orbit of a periodic point, which is a finite set, is called a \emph{periodic orbit}.  The set of periodic points of $f$ in $X$ is denoted by $P(f)$. The set of all limit points of $\cO(x)$ is called the \emph{omega-limit set} of $f$ at $x$, and written as $\omega_f(x) = \omega(x) $.  Omega-limit sets are closed  invariant sets.
 An element $x \in X$ is called a \emph{recurrent point} for $f$ if for some $n_k \nearrow \infty$, $f^{n_k}(x) \to x$, i.e. $x \in \o(x)$. The set of recurrent points of $f$ in $X$ is denoted by $\cR(f)$. An element $x \in X$ is called a \emph{nonwandering  point} for $f$ if for every open $U \ni x$, $\ \exists \ n \in \N$ such that $f^{n}(U) \cap U \neq \emptyset$. The set of all nonwandering  points of $f$ is denoted by $\Omega(f)$. The system $(X,f)$ is said to be \emph{non wandering} if $X = \Omega(f)$.

\bigskip

The system $(X,f)$ is said to be \emph{point transitive} if there is an $x_0 \in X$ such that $\overline{\mathcal{O}(x_0)} = X$. These points with dense orbits are called \emph{transitive points}.

The system is called \emph{topologically transitive} when  for every
pair of  nonempty, open sets $U,V \subset X$, there exists $n \in \N$ such that $f^n(U) \cap V \neq \emptyset$. Notice that $f$ is surjective, and so a nonempty, open $U$  implies $f^{-n}(U)$ is nonempty and open for every $n \in \N$.

These definitions of point transitivity and topological transitivity  are equivalent on all perfect, compact metric spaces.

Many times we represent these systems as $(X,x_0,f)$ or $(X,x_0)$ where $\overline{\mathcal{O}(x_0)} = X$. Such systems are then termed \emph{pointed systems or ambits}.

 The system is \emph{minimal} when every orbit is dense.

\bigskip


$(X,f)$ is called \emph{topologically mixing} if for every pair $V, W$ of nonempty open sets in $X$, there is a  $N>0$ such that $f^n(V) \cap W$ is nonempty for all $n \geq N$.

\bigskip

For $U, V \subseteq X$, let $N (U, V ) = \{n \in \N : f^n(U) \cap V \neq \emptyset\}$ be the \emph{hitting time set}. We say that

\begin{itemize}
	\item $(X,f)$ is transitive if for every pair of nonempty open sets $U, V \subseteq X$, $N (U, V )$ is nonempty.
%
	
	\item $(X,f)$ is mixing if for every pair of nonempty open sets $U, V \subseteq X,$ we have that $N (U, V )$ is cofinite.
	
\end{itemize}

\bigskip

An \emph{equivariant map}
$\pi : (X_1,f_1) \to (X_2,f_2)$ is a continuous map
$\pi : X_1 \to X_2$ such that $f_2 \circ \pi = \pi \circ f_1$.

In particular, the diagram

$$\begin{CD}
X_1  @>f_1>> X_1\\
@V{\pi}VV   @VV{\pi}V\\
X_2 @>f_2>> X_2
\end{CD}$$

  commutes.

  When $\pi$ is a homeomorphism we call it a \emph{conjugacy} and say that  $(X_1,f_1)$
and $(X_2,f_2)$ are \emph{conjugate}. When $\pi$ is surjective we call it a \emph{factor map} and say that $(X_2,f_2)$
is a \emph{factor} of $(X_1,f_1)$.

We note that the properties of entropy, transitivity and topologically mixing are preserved on taking factors.

\bigskip

We refer to \cite{ea, gh, v1} for more details on topological dynamics.

\bigskip

In \cite{hofer}, the concept of entropy has been extended by \textsc{Hofer} to non-compact Hausdorff spaces. According to \cite{hofer}, for a non-compact space $X$ and $T: X \rightarrow X$, the topological entropy $$h_{X}(T) = h_{X^*}(T^*),$$ where $X^*$ is a compactification of $X$ and $T^*$ is the extension of $T$ on $X^*$. We note that with this definition, it has been shown in \cite{hofer}:

(1) $h_X(T^k)=kh_X(T)$ for each positive integer $k$,

(2) If $Y \subset X$ is $T-invariant$, then $h_Y(T|_Y) \leq h_X(T)$,

(3) For non-compact spaces $X$ and $Y$ and systems $(X,T)$ and $(Y,S)$, if $\phi: X\rightarrow Y$ satisfies $\phi \circ T = S \circ \phi $ then $h_Y(S) \leq h_X(T)$.

\bigskip

It has been seen that pseudo-orbits, or more formally $\ep$-chains, are important tools for investigating properties
of discrete dynamical systems. They usually capture the recurrent
and mixing behaviors  depicted by the systems, even though they are highly metric dependent properties.

\smallskip

For $x,y \in X$, an \emph{$\ep$-chain (or $\ep$-
pseudo-orbit)} from $x$ to $y$ is a sequence $\{x = x_0, x_1, \ldots , x_n = y\}$ such
that $d(f(x_{i-1}), x_i) \leq \ep$ for $i = 1, \ldots , n$. The length of the $\ep$-chain $\{x_0, x_1, \ldots , x_n\}$ is said to be $n$.

A point $x \in X$ is \emph{chain recurrent} if for every $\ep > 0$,
there is an $\ep$-chain from $x$ to itself. $(X,f)$ is \emph{chain recurrent} if every point of $X$ is chain recurrent.

$(X,f)$ is \emph{chain transitive} if for every $x, y \in X$ and every $\ep > 0$, there is an $\ep$-chain from $x$ to $y$. $(X,f)$ is \emph{chain mixing} if for every $\ep > 0$ and  for any $x, y \in X$,  there is an $N > 0$ such that for all $n \geq N$, there is an $\ep$-chain from $x$ to $y$ of length exactly $n$.

\bt \cite{ea} \label{cto} For $(X,f)$,  $\o(x)$  is chain transitive for all $x \in X$. Further, if $(X,f)$ is chain transitive then it can be embedded in a larger system where it is an omega limit set. \et

\bt \cite{aw} \label{dicho}If $(X,f)$ is chain transitive, then either $(X,f)$ is chain mixing or $(X,f)$ factors onto a non-trivial periodic orbit. \et

We suggest to the enthusiastic readers to look into \cite{ea, aw, rw} and the references therein for many interesting  chain properties.

\subsection{Symbolic Dynamics}

Shift spaces are built on a finite set $\mathcal{A}$ of symbols which we call the \emph{alphabet}. Elements of $\mathcal{A}$ are called \emph{letters}.

We define the \emph{full  $\mathcal{A}$  shift} as the collection of all bi-infinite sequences or bi-infinite sequences of symbols from $\mathcal{A}$. It is denoted by

\begin{align*} \mathcal{A}^{\mathbb{Z}} = \{ x\ =\ ...x_{-1}.x_0x_1...:\ x_i\ \in\ \mathcal{A}\ \forall\ i \in \mathbb{Z} \}.\end{align*}

The product topology on $\A^{\Z}$ is metrizable and a compatible metric defined on it can be given as:

\begin{equation} \label{metric}
d(x,y) = \inf \left\lbrace  \frac{1}{2^m}  : x_n = y_n \ \text{for} \ |n| < m  \right\rbrace,
\end{equation}

for any two sequences $x = ...x_{-1}.x_0x_1... $ and $y = ...y_{-1}.y_0y_1... \in \A^{\Z}$.

\bigskip

The $shift\ map \ \sigma$ on the full shift $\mathcal{A}^{\mathbb{Z}}$ maps a point $x$ to the point $\sigma (x)$ whose ith coordinate is \begin{align*} (\sigma(x))_i\ =\ x_{i+1}. \end{align*}

A \emph{shift space} is a closed, invariant (i.e. $\sigma(X) \subseteq X$) set  $X \subseteq \mathcal{A} ^ {\mathbb{Z}}$.

 Observe that for $x, y \in X$,

$$\exists \ m>0 \ \text{for which} \ d(x, y) < 2^{-m} \Leftrightarrow \ \exists \ k>0 \ \text{for which} \ x_{[-k,k]} = y_{[-k,k]}.$$

\bigskip

For $n \in \N$,  (a nonempty) $w \in \A^n$ is a word of length $n$, and we write $|w| = n$. If the word $w$ is a part of the word $v$ then we say that $w$ is a subword of $v$  and we write $w \sqsubset v$.

Similarly for any $x \in \A^{\Z}$, we write $w \sqsubset x$ if $w$ appears in $x$ as a  block, i.e. $w = x_{[k, k+n]} = x_kx_{k+1} \ldots x_{k+n} $ and for $m \in  \N$, we write the concatenation $w^m = \underbrace{w \ldots w}_{m \text{-times}}$. The collection of all nonempty words in $\A^{\Z}$ is  $\A^* = \bigcup \limits_{n \in \N} \A^n$.

Let $\cL(X) \subset \A^*$ be the language of  shift space $X$ i.e. the set of all nonempty words appearing in any $x \in X$. Let $\mathcal{F} \subset \A^*$ be the set of blocks that never appear in any $x \in X$. Usually we can write the shift space $X$ as $ X_{\mathcal{F}} $ for some collection $\mathcal{F}$ of forbidden blocks over $\mathcal{A}$, i.e.  $ \mathcal{F} \subseteq \mathcal{L}(X)^c  $. Every shift space can also be defined by its language $X = X_{\mathcal{F}} = X_{{\mathcal{L}(X)}^c} $. Notice that for $X = \A^{\Z}$ we have $\cF = \emptyset$.

\bigskip

For the shift space $X \subseteq \A^{\Z}$, the system $(X,\s)$ is called a \emph{subshift}.

\bigskip

$( X_{\mathcal{F}}, \s) $ is called a \emph{subshift of finite type (SFT)} if the list of forbidden words $\mathcal{F}$  can be taken to be finite. If $M + 1$ is the length of the longest forbidden word, then this SFT is an $M-$step
SFT. Thus an $M-$step SFT has the property that if
$uv, vw \in \mathcal{L}(X)$  and $|v| \geq M+1$, then $uvw \in \mathcal{L}(X)$ as well.

Every SFT is topologically conjugate to an
\emph{edge SFT} $X_A$, presented by some square
nonnegative matrix $A$. Here $A$ is viewed as
the adjacency matrix of some directed graph $G$,
whose edge set is the alphabet $\A$ of the SFT. $X_A \subset \A^{\Z}$  is the
space of bi-infinite sequences  corresponding to walks through the graph $G$. Here for every $i$, the terminal vertex of $x_i$ equals
the initial vertex of $x_{i+1}$.  Thus, a SFT can also be denoted as $X_G$ or $X_A$ where $G$ is the associated graph and $A$ is the transition matrix.

\bigskip

A non-negative matrix $A = (a_{ij})$ is called \emph{irreducible} if for every $i, j$ there is $k \in \N$ such that $A^k_{ij} > 0$, and is called \emph{aperiodic} if there is $k \in \N$ such that $(A^k)_{ij} > 0$, for every $i, j$.

A subshift $(X,\s)$ is transitive if for every pair of words $u, v \in \cL(X)$, there exists a word $w \in \cL(X)$ with $uwv \in \cL(X)$. $(X,\s)$ is mixing if there exists an $N \in \N$ such that such a $w$ can be choosen with $|w|=n$ for all $n \geq N$.

A SFT $(X_A, \s)$ is transitive if and only if the transition matrix $A$ is irreducible, and is mixing if and only if the transition matrix $A$ is aperiodic.

SFTs can also be viewed as \emph{vertex shifts}, where the vertices are labelled by elements in $\A^*$ and a directed edge connects two vertices if the concatenation of the labels of these vertices is a permissible word in the language. The transition matrix here is a Boolean matrix.
\bigskip

Here codes play an important role. The most important codes for us are those that do not change with time i.e. codes which intertwine with the shift $(\sigma \ \circ \ \phi\ =\ \phi \ \circ \ \sigma)$. A common example of such codes is \emph{the sliding block code} which we define as follows: Let X be a shift space over $\mathcal{A}$. Let $\mathcal{B}$ be another alphabet and \begin{align*} \bar{\phi} : \mathrm{B}_{m+n+1} (X)\ \rightarrow\ \mathcal{B}\end{align*} be a map called \emph{(m+n+1)-block map} or simply \emph{block map}, where $\mathrm{B}_{m+n+1} (X)$ is the set of all $(m+n+1)$-blocks in $\cL(X)$. We define a map \begin{align*} \phi : X \rightarrow \mathcal{B}^{\mathbb{Z}}\end{align*} given by $y = \phi(x)$ where $y_i\ =\ \bar{\phi} (x_{[i-m,i+n]})$. The map $\phi$ is called the \emph{sliding block code} with \emph{memory} m and \emph{anticipation} n induced by $\bar{\phi}$. \par

\bt (Curtis-Hedlund-Lyndon)  Let $(X, \s_X)$ and $(Y, \s_Y)$ be subshifts over finite alphabets $\mathcal{A}$ and $\mathcal{B}$ respectively. A continuous map $\phi\ :\ X\ \rightarrow\ Y$ commutes with the shift
	(i.e., $\phi \ \circ \ \sigma_X\ =\ \sigma_Y \ \circ \ \phi$) if and only if $\phi$ is a sliding block code.\et

If a sliding block code $\phi\ :\ X\ \rightarrow\ Y$ is onto, it is called a \emph{factor code} from X onto Y and Y is called a \emph{factor} of X. If $\phi$ is one-to-one, then it is called \emph{embedding} of X into Y. $\phi$ is called a \emph{conjugacy} from X to Y, if it is invertible. The shift spaces in this case are called $conjugate$ and we write X $\equiv$ Y. Conjugacies carry $n-$periodic points to $n-$periodic points and in general preserve the dynamical structure.

\bigskip

The topological entropy of a subshift $X$ is given as 
 $$h_{X} = \lim \limits_{n \to \infty} \frac{\log | \mathrm{B}_{n}|}{n},$$
  where $|\mathrm{B}_{n}|$ denotes the number of  words in $\mathcal{L}(X)$ of length $n$. It is known that the topological entropy of an irreducible SFT $X_{A}$ equals $\log \lambda$ where $\lambda$ is the \emph{Perron eigenvalue} of $A$.

\bigskip

We refer to \cite{marcus} for more details.

\bigskip

There is an interesting illustration of Theorem \ref{cto} for SFTs.

\bt \cite{bgkr} Let $\Lambda \subset X_{\cF}$ be an invariant
and closed subset. Then there is a point $x \in X_{\cF}$ such that $\Lambda = \o(x)$ if and only if $\Lambda$ is chain transitive. \et

\br \label{ct=t} One of the consequences of the above observation, also observed independently in \cite{rw} is - \emph{in the case of  SFTs, chain transitivity is equivalent to transitivity}. \er

Lastly, we recall some characterization of the language of a subshift from \cite{eja} and build on it. We skip the trivial proofs since the arguments are similar to the ones given in \cite{eja}. Recall that as $v$
varies over $\cL(X)$, the \emph{cylinder sets}
$$[v]\ = \ \{ x \in X : x_{[-k, -k + |v|]} = v, k \in \N \}\hspace{2cm}$$
comprise a bases of clopen sets  on $X$.

\begin{itemize}
	\item[(a)] $(X, \s)$ is  transitive(irreducible) if and only if for all $v \in \cL(X)$ and  all $w \in \cL(X)$, there
	exists  $a \in \cL(X)$ such that $vaw \in \cL(X)$.
	
	\item[(b)] $(X, \s)$ is minimal  if and only if whenever $v \in \cL(X)$ then $v \sqsubset x$ for all $x \in X$.
	
	\item[(c)] $(X, \s)$ is  mixing if and only if whenever $v, w \in \cL(X)$ there
	exists $N \in \N$ such that for all $k \in \N$ there exists $a_k \in \cL(X)$ with $|a_k| = N + k$ such that $va_kw \in \cL(X)$.
	
	\item[(d)] $(X, \s)$ is non wandering if and only if for all $v \in \cL(X)$ there	exists  $a \in \cL(X)$, such that $vav \in \cL(X)$.
	
	\item[(e)] $(X, \s)$ is chain recurrent if and only if for all $v \in \cL(X)$, there	exists  $a_1, a_2, \ldots a_n \in \A$ and  $v_1, v_2, \ldots v_{n-1} \in \cL(X)$ with $|v_1| = |v_2| = \ldots = |v_{n-1}|$, such that $va_1v_1, v_1a_2v_2, \ \ldots, \ v_{n-1}a_nv \in \cL(X)$.
	
	\item[(f)] $(X, \s)$ is chain transitive if and only if for all $v \in \cL(X)$ and  all $w \in \cL(X)$, there	exists  $a_1, a_2, \ldots a_n \in \A$ and  $v_1, v_2, \ldots v_{n-1} \in \cL(X)$  with $|v_1| = |v_2| = \ldots = |v_{n-1}|$, such that $va_1v_1, v_1a_2v_2, \ \ldots, \ v_{n-1}a_nw \in \cL(X)$.
	
	\item[(g)] $(X, \s)$ is chain mixing if and only if for all $v \in \cL(X)$ and  all $w \in \cL(X)$,  there
	exists $N \in \N$ such that for all $k \in \N$ there exists $a_{m_j} \in \A$ for $j=1, \ldots, N+k$   and  $v_1, v_2, \ldots v_{N+k-1} \in \cL(X)$  with $|v_1| = |v_2| = \ldots = |v_{N+k-1}|$, such that $va_{m_1}v_1, v_1a_{m_2}v_2, \ \ldots, \ v_{N+k-1}a_{m_{N+k}}w \in \cL(X)$.
	
\end{itemize}

\subsection{Some Rudiments in Metric Spaces}

For a metric space $(X,d)$, we  denote as
$2^X$ $-$  the space of all nonempty closed subsets of $X$, and $\K(X)$ $-$ the space of all compact subsets of $X$, endowed with the Hausdorff topology.  We note that usually $\K(X) \subseteq 2^X$ but for compact $X$, $2^X = \K(X)$.

This has a natural induced metric.

\smallskip

Given a point $p \in X$ and a closed set $A \subseteq X$, recall
$d(p, A) = \inf \limits_{a \in A} d(p, a).$

\smallskip

On $2^X$ we define \emph{the Hausdorff metric}:  For $A, B \in 2^X$

\begin{equation}\label{1}
d_H(A,B) \quad = \quad \max \{\sup \limits_{a \in A} d(a, B),  \sup \limits_{b \in B} d(b, A)\}
\end{equation}

We note that $d_H$ is a psuedo-metric on $2^X$ and a metric on $\K(X)$.

\bigskip

For $\ep > 0$, let $A_{\ep} = \{y \in X: d(y,a) < \ep, \ for \ some\  a\ \in\ A\}$ be the $\ep-$neighbourhood of $A$. Thus, $d_H(A,B) < \ep $ if and only
if each set is in the open $\ep $ neighborhood of the other i.e. $A \subset B_\ep$ \ and \ $ B \subset A_\ep$,
or, equivalently, each point of $A$ is within $\ep$ of a point in $B$ and vice-versa.

\bigskip

When $X$ is compact, we occasionally use an equivalent topology on $2^X$. Define for
any collection $\{ U_i : 1 \leq i \leq n \}$ of  open and nonempty subsets of $X$,

\begin{equation}\label{1.1a}
\langle U_1, U_2, \ldots U_n \rangle =  \{E \in 2^X :E \subseteq \bigcup \limits_{i=1}^n
U_{i}, \ E \bigcap U_{i} \neq \phi,  \textrm{ } 1 \leq i \leq n \}
\end{equation}
The topology on $2^X$, generated by such collection as basis, is
known as the {\emph{Vietoris topology}}.

\bigskip

If $\{ A_n \}$ is a sequence of closed sets in a $X$ then
\begin{equation*}
\begin{split}
\overline{ \bigcup_n \{ A_n \}} \quad = \quad \bigcup_n \{ A_n \} \ \cup \ \limsup_n \{ A_n \}, \hspace{1cm}\\
\mbox{where} \qquad \ \limsup_n \{ A_n \} \quad = \quad \ \bigcap_k \overline{\bigcup_{n \geq k} \{ A_n \}}.
\end{split}
\end{equation*}

We recall,

\bl  For a  metric space $(X,d)$, and $ \{ A_n \}$  a sequence in $(2^X, d_H)$,
	
	(a) If $\{ A_n \}$ converges to $A$ in $\K(X)$ then $\overline{ \bigcup_n \{ A_n \}}$ is compact.
	
	(b) If $\overline{ \bigcup_n \{ A_n \}}$ is compact and $\{ A_n \}$ is Cauchy then $A_n$ converges to  $ \limsup \{ A_n \}$.
	
	(c) If $X$ is complete then $2^X$ is complete.
	
	(d) If $X$ is compact then $2^X$ is compact.
	
	(e) If $X$ is separable then so is $2^X$.
	
	(f) If $i_X : X \to 2^X$ is given by $x \mapsto \{ x \}$, then $i_X$ is an  isometric inclusion.
	
	(g) For $f : X_1 \to X_2$ a continuous map of metric spaces, there is induced the map $f_* : \K({X_1}) \to \K({X_2})$  defined by $A \mapsto f(A)$. If $f$ is
	uniformly continuous or continuous, then the map $f_*$ is also uniformly continuous or continuous, respectively.
	
%
	(h) The set of all finite subsets of $X$ is dense in $\K(X)$.
	
\el

These results also hold when $2^X$ is given the Vietoris topology.

\bigskip

We refer  \cite{aan, in, mi}  for more details.

\bigskip

The \emph{Gromov-Hausdorff metric} furthers the idea of the Hausdorff metric.
Given two compact metric spaces $X$ and $Y$, we define

\begin{equation} \label{2}
d_{GH}(X, Y) = \inf \limits_{f,g} d_H(f(X), g(Y))
\end{equation}

where $f(X), g(Y)$ denote an isometric embedding of $X, Y$ into some metric space $Z$  and the infimum is taken
over all such possible  embeddings.

\bl \cite{gr} The following hold with respect to the Gromov-Hausdorff metric:

\be
\item If $X,Y$ are compact metric spaces, the $d_{GH}(X,Y) < \infty$.

\item If $X,Y$ are not compact then it is possible that $d_{GH}(X,Y)=0$, without $X,Y$ being isometric. For example $[0,1], \mathbb{Q} \cap [0,1]$.

\item Metric spaces  $X_i \to X$ if and only if for every $\ep > 0$, there exists $\ep' \geq \ep$ such that every $\ep'-$net in $X$ is a limit of $\ep-$nets in $X_i$.

\item Compact $X$ and $Y$ are isometric if and only if $d_{GH}(X, Y ) = 0$.

\ee

\el

We denote the set of  all isometric compact metric spaces endowed with the Gromov-Hausdorff metric  as  $\M$. Then the following is known about the metric space $(\mathcal{M}, d_{GH})$.

\be

\item $(\mathcal{M}, d_{GH})$ is separable and complete.

\item $(\mathcal{M}, d_{GH})$ is not locally compact or compact.

\item The set of all finite metric spaces is  dense in $(\mathcal{M}, d_{GH})$.

\ee

Such a $\mathcal{M}$  is called a \emph{universal metric space}.

\bigskip

We refer to \cite{bbi, gr} for more details.

\subsection{Geodesics and Polygons in the Hyperbolic Plane}

A hyperbolic space is a space that has a constant negative sectional curvature. We work in dimension $2$ and call it a hyperbolic 2-space or a hyperbolic plane. We will use two models of hyperbolic plane, namely, the Poincar\'e half plane model which we denote as $\mathbb{H}$ and the Poincar\'e disc model denoted as $\mathbb{D}$. Both of them model the same geometry in the sense that they can be related by an isometric transformation that preserves all the geometrical properties. We refer to \cite{anderson, beardon, series} for more details.\par

The underlying space of the Poincar\'e half plane model is the upper half-plane $\mathbb{H}$ in the complex plane $\mathbb{C}$, defined to be
$$\mathbb{H}\ = \{z\  \in\ \mathbb{C}\  |\ Im(z) > 0\}.$$ We use the usual notion of point and angle that $\mathbb{H}$ inherits from $\mathbb{C}$.
The metric on $\mathbb{H}$ is defined by $$ds^2 = \dfrac{dx^2 + dy^2}{y^2}.$$

\bigskip

The Poincar\'e disk model is described by $$\mathbb{D} =  \{z \ \in\ \mathbb{C}\ : \ |z| \ < \ 1\}.$$ The metric on $\mathbb{D}$ is defined by $$ds^2= \dfrac{4(dx^2+dy^2)}{(1-(x^2+y^2))^2}.$$ 

\smallskip

 Once we  have the notion of distance on a space, we define the respective \emph{straight lines} as \emph{geodesics}. They are the locally distance minimising curves of the space. Under the metric imposed on $\mathbb{H}$, we get the geodesics to be the euclidean lines perpendicular to real axis and the euclidean semicircles which are orthogonal to the real axis. We note that the real axis along with the point at infinity gives the boundary $\partial \mathbb{H}$ of $ \mathbb{H}$.  In case of $\mathbb{D}$ with the above defined metric, we get the geodesics to be the euclidean lines passing through the centre of the disc and the euclidean circles orthogonal to $\partial \mathbb{D}$, the boundary of $ \mathbb{D}$.

\bigskip

A subset $A$ of the hyperbolic plane is \emph{convex} if for each pair of distinct points $x$ and $y$ in $A$, the closed line segment $l_{xy}$ joining $x$ to $y$ is contained in $A$. Hyperbolic lines, hyperbolic rays, and hyperbolic segments are convex. Given a hyperbolic line $l$, the complement of $l$ in the hyperbolic plane has two components, which are the two \emph{open half-planes determined by l}. A \emph{closed half-plane determined by l} is the union of $l$ with one of the two open half-planes determined by $l$. We refer to $l$ as the \emph{bounding line} for the half-planes it determines. Open half-planes and closed half-planes in $\mathbb{H}$ are convex.

Let $\mathsf{H}=\{H_{\alpha} \}_{\alpha \in \Lambda}$ be a collection of half-planes in the hyperbolic plane, and for each $\alpha \in \Lambda$, let $l_{\alpha}$ be the bounding line for $H_{\alpha}$. The collection $\mathsf{H}$ is called \emph{locally finite} if for each point $z$ in the hyperbolic plane, there exists some $\epsilon > 0$ so that only finitely many bounding lines $l_{\alpha}$ of the half-planes in $\mathsf{H}$ intersect the open hyperbolic disc $U_{\epsilon}(z)$ where $U_{\epsilon}(z)= \{w \in \mathbb{H} : d_{\mathbb{H}}(z,w)< \epsilon\}.$

  A \emph{hyperbolic polygon} is a closed convex set in the hyperbolic plane that can be expressed as the intersection of a locally finite collection of closed half-planes. Under this definition, there are some subsets of the hyperbolic plane that satisfy this criteria, but we do not want them to be considered as hyperbolic polygons. For example, a hyperbolic line $l$ is a hyperbolic polygon, because it is a closed convex set in the hyperbolic plane that can be expressed as the intersection of the two closed half-planes determined by $l$. A hyperbolic polygon is \emph{nondegenerate} if it has nonempty interior else it is called \emph{degenerate}. We will work only with nondegenerate polygons here.
  
  \bigskip

Let $P$ be a hyperbolic polygon and let $l$ be a hyperbolic line so that $P$ intersects $l$ and so that $P$ is contained in a closed half-plane determined by $l$. If the intersection $P \cap l$ is a point, we say that this point is a \emph{vertex} of $P$. The other possibilities are that the intersection $P \cap l$ is either a closed hyperbolic line segment, a closed hyperbolic ray, or all of $l$. We call this intersection \emph{side} of $P$. Let $P$ be a hyperbolic polygon, and let $v$ be a vertex of $P$ that is the intersection of two sides $s_1$ and $s_2$ of $P$. Let $l_k$ be the hyperbolic line containing $s_k$. The union $l_1 \cup l_2$ divides the hyperbolic plane into four components, one of which contains $P$. The \emph{interior angle} of $P$ at $v$ is the angle between
$l_1$ and $l_2$, measured in the component of the complement of $l_1 \cup l_2$ containing $P$. A hyperbolic polygon $P$ in the hyperbolic plane has an \emph{ideal vertex} at $v$ if there are two adjacent sides of $P$ that are either closed hyperbolic rays or hyperbolic lines  that share $v$ as an endpoint at infinity.

\bigskip

A finite-sided polygon $P$ in the hyperbolic plane is called \emph{reasonable} if $P$ does not contain an open half-plane. A \emph{hyperbolic n-gon} is a reasonable hyperbolic polygon with $n$ sides. A \emph{compact polygon} is a hyperbolic polygon whose all vertices are in the hyperbolic plane. A compact hyperbolic n-gon is \emph{regular} if its sides have equal length and if its interior angles are equal. For each $n \geq 3$, an \emph{ideal\ n-gon} is a reasonable hyperbolic polygon $P$ that has $n$ sides and $n$ vertices. Thus, an ideal polygon is a hyperbolic polygon whose all vertices are ideal points (i.e. lying on the boundary of the hyperbolic plane). The hyperbolic polygons with vertices lying both inside the hyperbolic plane and on its boundary are called \emph{semi-ideal polygons}. An angle at a vertex of a polygon in $\mathbb{D}$ is called \emph{rational} if it is of the form $\pi / n$ where $n \in \mathbb{N}$ and $n > 1$. The corresponding vertex is called a \emph{rational vertex}. Thus, a compact polygon is labeled \emph{rational} if all its vertices are rational and a semi-ideal polygon is called \emph{rational} if all its non-zero vertex angles are rational. We note that  the angle at vertices that are ideal points is zero and so by definition ideal polygons are vacuous rational polygons. 

\bigskip

An ideal polygon has infinite perimeter and finite area from Gauss-Bonnet formula. In particular, an ideal $(k+2)$-sided polygon has an area $k\pi$ and thus is the largest possible polygon in hyperbolic plane. The compact polygons have finite perimeter and area strictly less than $k\pi, \ k\in \N$. The semi-ideal ones have infinite perimeter and area less than or equal to $k\pi, \ k \in \N$. The polygons in the hyperbolic plane enjoy a very special feature, namely, the similar polygons are congruent. In particular, all ideal n-gons are congruent to each other. This feature allows us to work in a simplified situation of a symmetrically placed n-gon.

\bigskip

We refer to \cite{anderson,loring} for more details.

\bigskip

\subsection{Tessellating the Hyperbolic Plane and the Katok-Zemlyakov unfolding method}

We can consider any one of $\mathbb{H}\ or\ \mathbb{D}$ as the model for the hyperbolic plane. A \emph{tessellation} of $\mathbb{D}$ is a subdivision of $\mathbb{D}$ into polygonal tiles $\Pi_i,\ i \in \Lambda$ satisfying the following conditions:

 (1) $\forall\ z \in \mathbb{D},\ \exists\ i \in \Lambda $ such that $z \in \Pi_i$,
 
  (2) $\forall\ i \neq j,\ \Pi_i \cap \Pi_j$ is either empty, or a single vertex common to both, or an entire common edge,
  
   (3) $\forall\ i \neq j,\ \exists $ an isometry $f_{i,j}$ of hyperbolic plane  such that $f_{i,j}(\Pi_i) = \Pi_j$.

   Informally speaking, a collection of tiles tessellate the hyperbolic plane if they cover the plane, don't overlap, and are of same shape and size.


Let $\Pi$ be an ideal polygon in the hyperbolic plane. Then, we can reflect it across each one of its sides and the same procedure can be applied to the reflections and so on. The collection of all such ideal polygons obtained, along with $\Pi$ gives us a tessellation of $\mathbb{D}$. The figure below shows one such example, where we start with an ideal triangle $\Pi$. Its sides are labeled as $1, 2, 3$ in counter-clockwise sense. On reflection about a side $i$, the labels change to $1^i,2^i,3^i$. This labelling proceeds in same way for further reflections. Similar tessellation can be obtained for the rational compact polygons and for the semi-ideal polygons with the vertex angles either $0$ or rational.

\begin{figure}[h!]
	\centering
	\includegraphics[width=8cm,height=7.5cm]{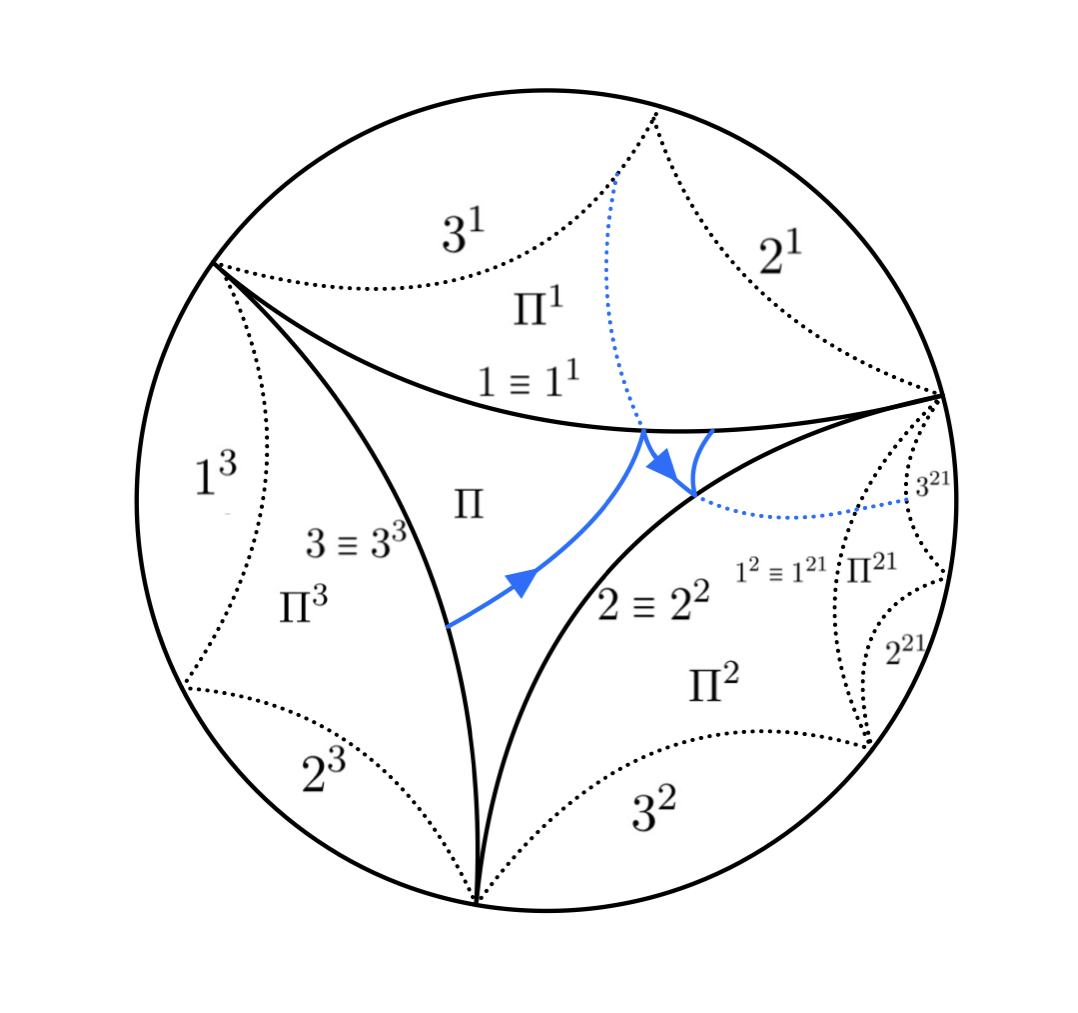}
	\caption{Tessellation of $\mathbb{D}$ and unfolding\ of\ a\ billiard\ trajectory}
\end{figure}

The unfolding technique that we discuss ahead was formally introduced to the domain of billiard dynamics by 
  \textsc{A.B.Katok} and \textsc{A.N.Zemlyakov} \cite{katok}. This method converts a polygon on a plane into a surface on which the billiard trajectories appear as geodesics.
Under this technique a billiard trajectory in a polygon $\Pi$ can be \emph{unfolded} in an intuitive procedure as follows: Instead of reflecting the trajectory in a side of $\Pi$, we reflect $\Pi$ itself in that side which gives a copy of the reflected ray in the new polygon. The join of this new directed line segment with the incident ray in $\Pi$ lies on a line. We say that the billiard trajectory has been \emph{unfolded} at this hit point. When we apply this procedure to the whole billiard trajectory, it gives us straightened version of the billiard trajectory, which we call as \emph{unfolded billiard trajectory}.

\bigskip

More details on tessellations of polygons in the hyperbolic plane can be found in the elaborated work of S.Katok in  \cite{skatok}, wherein it has been described via the Fuchsian groups and the fundamental regions. The tessellations are described via the action of a discrete and properly discontinuous group on a fundamental region. Here, we avoid this nomenclature as our primary object of concern is the polygon (which is already fixed to start with) and the billiard dynamics happening on it.

\bigskip

We refer to \cite{beardon,chernov, tabachnikov} for more details.

\subsection{Billiards in the Hyperbolic Plane} 

We consider the Poincar\'e disc model $\mathbb{D}$ of the hyperbolic plane.

\bigskip

We carry out our investigations with a class of hyperbolic polygons  which we call as \emph{ideal polygons, compact rational polygons and semi-ideal rational polygons}. The \emph{ideal polygons} are the ones with all the vertices on $\partial \mathbb{D}$ and the \emph{compact rational polygons} are the ones for which all the vertices lie in $\mathbb{D}$ and have angles that divide $\pi$ into integer parts. The polygons that lie in `semi-ideal rational' class are the ones for which either the vertices lie on $\partial \mathbb{D}$ or the vertex angles divide $\pi$ into integer parts. 

\bigskip

 We follow the general construction as defined in  \cite{ullmo} for the ideal polygons and in \cite{ullmo2} for compact rational polygons.
 
 \bigskip
 
 Let $\Pi$ be a k-sided polygon in $\mathbb{D}$. A billiard trajectory in the polygon $\Pi$ is a directed geodesic flow between each pair of consecutive specular bounces of the boundary(not containing the vertices) of $\Pi$. We have a simple choice for the coordinate system for the directed geodesic arcs. We parameterize the boundary of the Poincar\'e disc $\mathbb{D}$ using the azimuthal angle by considering it as a subset of $\mathbb{C}$. Thus, we can represent a directed geodesic by the pair $(\theta , \phi )$, where $\theta , \phi$ are the intercepts made by a directed geodesic on $\partial \mathbb{D}$ with the direction being from $\theta$ to $\phi$. In this setting we have a natural metric on $\partial \mathbb{D}$ given by
 
\begin{equation} \label{dbg}
  d_{\partial \mathbb{D}} (\phi_1 , \phi_2) = |\phi_1 - \phi_2 |.\end{equation}

A billiard trajectory is a curve that is parameterized by the arc-length and consists of the geodesic arcs which are reflected by the walls of the polygon $\Pi$. Therefore, a trajectory can be expressed as \begin{align*} \gamma = \{(\theta_n , \phi_n)\label{key}_{n \in \mathbb{Z}}\} \end{align*} where $$(\theta_n , \phi_n) = T(\theta_{n-1} , \phi_{n-1}).$$

 We do not consider billiard trajectories starting or ending in vertices of $\Pi$. 
 
 \bigskip
 
 The billiard trajectories here follow the grammar rules defined and discussed by  \textsc{Marie-Joya Giannoni, Dennis Ullmo} in \cite{ullmo,ullmo2}. Based on that \textsc{Simon Castle, Norbert Peyerimhoff, Karl Friedrich Siburg} \cite{castle} prove the rules for ideal polygons.

\bt
\cite{castle} \label{castle result} Let $\Pi  \subset \mathbb{D}$ be  an  ideal  polygon   with  counter-clockwise enumeration  $1,...,k$. An equivalence class  $[...a_{-1}.a_0 a_1...]$ denoted $(a_j)$ with $...a_{-1}.a_0 a_1...  \in\ \{1,...,k\}^{\mathbb{Z}}$  is  in  $S(\Pi)$ if and  only  if

(1) $(a_j)$ does  not  contain immediate   repetitions,  i.e., $a_j  \neq  \ a_{j+1} \ \forall \ j \ \in\ \mathbb{Z}$ and

(2) $(a_j)$ does  not  contain  an  infinitely  repeated  sequence  of  labels  of  two  adjacent  sides.

Moreover,  every  equivalence class of pointed billiard  sequences corresponds  to  one and  only one billiard trajectory. \et

We refer \cite{castle,ullmo,ullmo2} for further details on the geometric aspects of billiards in hyperbolic plane.

\section{Pointed Geodesics and Billiards in Hyperbolic Polygons}

\subsection{Pointed Geodesics and Billiards}
\bd Let $\gamma = (\theta_n,\phi_n)_{n \in \mathbb{Z}}$ be a \emph{billiard trajectory} in a polygon $\Pi$ in $\mathbb{D}$. For a fixed $n \in \mathbb{Z}$, we will call $(\theta_n,\phi_n)$ as a \emph{base arc} of the trajectory $\gamma$.\ed 

We note that base arcs are compact subsets of $\mathbb{D}$.

\smallskip

A base arc uniquely determines the billiard trajectory under the restrictions imposed by the specular reflection rule.

\bd For the base arc $(\theta, \phi)$ defining $\gamma $, we call $(\gamma,(\theta,\phi))$ a \emph{pointed geodesic}. \ed 

Thus a pointed geodesic $(\gamma,(\theta,\phi))$ is identified with the element $$\ldots (T^{-1}(\theta,\phi)).(\theta,\phi)(T(\theta,\phi)) \ldots \in \K(\mathbb{D})^{\Z}$$

 by clearly establishing the position of the base arc $(\theta,\phi)$. A natural way of encoding a pointed geodesic is to seize the order in which it hits the sides of $\Pi$, starting from the side hit by the base arc and then reading the past and future hits of the trajectory and pointing out the symbol corresponding to the base arc. If we label the sides of $\Pi$ anti-clockwise from 1 to k, then every pointed geodesic produces a  bi-infinite sequence $...a_{-1}.a_0a_1...$ with $a_j\ \in\ \{1,...,k \}$.
 
 \bd  Define \begin{align*} \mathbb{G} = \mathbb{G}_{\Pi}= \{ \big(\gamma, (\theta,\phi) \big) : \gamma  = \big( T^n(\theta,\phi) \big)_{n \in \mathbb{Z}} \} \end{align*} as the \emph{space of all pointed geodesics} on $\Pi$.
 
  $\mathbb{G} \subseteq \K({\mathbb{D}})$ and so    ${\mathbb{G}}$ can be equipped  with the natural Hausdorff metric $d_H$, and so is endowed with the Hausdorff topology. \ed

 Here, $ T^0 (\theta, \phi)$ is simply written as $(\theta, \phi)$. 

\bigskip

We define a function $ d_{\mathbb{G}} : \mathbb{G} \times \mathbb{G} \to \mathbb{R}$ as follows:
 \begin{equation} \label{dg}
  d_{\mathbb{G}} \Big(\big(\gamma, (\theta,\phi) \big),\big(\gamma', (\theta',\phi') \big)\Big) = \max\{ d_{\partial\mathbb{D}}(\theta,\theta'),d_{\partial\mathbb{D}}(\phi,\phi')\}. \end{equation} 
  
  where $d_{\partial\mathbb{D}}$ is as in Equation (\ref{dbg}).

\bpr Let $\mathbb{G}$ be the space of pointed geodesics on a polygon $\Pi$ in $\mathbb{D}$, then $d_{\mathbb{G}}$ defines a metric on ${\mathbb{G}}$. \epr
\bo  Clearly, $d_{\mathbb{G}}$ is non-negative.
If $d_{\mathbb{G}} \Big(\big(\gamma, (\theta,\phi) \big),\big(\gamma', (\theta',\phi') \big)\Big) = 0,$ then
$ d_{\partial\mathbb{D}}(\theta,\theta') = 0 , d_{\partial\mathbb{D}}(\phi,\phi') = 0.$ Therefore,
$ \theta = \theta' , \phi = \phi' $ implying $ (\theta,\phi) = (\theta',\phi')$.
If base arcs of two pointed geodesics are same then the corresponding trajectories are also same because of the dynamics provided by the bounce map. Therefore, $$\Big(\big(\gamma, (\theta,\phi) \big)\Big) = \Big(\big(\gamma', (\theta',\phi') \big)\Big).$$ The symmetry and triangle inequality for $d_{\mathbb{G}}$ follows from the respective properties of $d_{\partial\mathbb{D}}$(the boundary of the Poincar\'e disc).
Therefore, $d_{\mathbb{G}}$ is a metric on ${\mathbb{G}}$.\eo

 We show that the Hausdorff topology on $\mathbb{G}$ is same as the topology on ${\mathbb{G}}$ given by $d_{\mathbb{G}}$. Recall that $d_H$ on $\mathbb{G}$ can be given as follows:

\begin{equation} d_H \Big(\big(\gamma, (\theta,\phi) \big),\big(\gamma', (\theta',\phi') \big)\Big)\ := \ d_H((\theta,\phi),(\theta',\phi'))\ = 
	\max\lbrace \displaystyle \sup_{Q \in (\theta,\phi)} d\big(Q,(\theta',\phi')\big), \displaystyle \sup_{Q \in (\theta',\phi')} d\big(Q,(\theta,\phi)\big)\rbrace
 \end{equation}

Note that the above definition works because $\gamma$ is uniquely determined by a base arc and if $\gamma \neq \gamma'$ then $d_H(\gamma,\gamma') > 0$. Thus, this notion of distance between two pointed geodesics is just the Hausdorff distance between the corresponding base arcs. We denote the space of all base arcs on a polygon $\Pi$ that are associated with billiard trajectories by $\mathbb{B}(\Pi)$ or simply $\mathbb{B}$, when the context is clear. Note that $\mathbb{B} \subset \K(\mathbb{D})$ and is a bounded subset of $\mathbb{D}$. Thus, in this case, the Vietoris topology and Hausdorff topology are equivalent on $\mathbb{B}$.      We, thereby get a natural isometry between $(\mathbb{G},d_H)$ and $(\mathbb{B},d_H)$ for a polygon $\Pi$, giving a one-to-one correspondence between the Vietoris topology on $\mathbb{B}$ and the topology generated by $d_H$ on $\mathbb{G}$. Under the same pretence, we also have the $d_{\mathbb{G}}$ metric on $\mathbb{B}$ and the natural isometry between $(\mathbb{G},d_{\mathbb{G}})$ and $(\mathbb{B},d_{\mathbb{G}})$.

\bt Let $\mathbb{G}$ be the space of pointed geodesics on a polygon $\Pi$ in $\mathbb{D}$, then  $d_{\mathbb{G}}$ and $d_H$ generate the same topology on $\mathbb{G}$. \et
\bo With the above discussion, it is sufficient to prove that $d_{\mathbb{G}}$ and $d_H$ generate the same topology on $\mathbb{B}$. The topology on $\mathbb{B}$ given by the metric $d_H$ is the induced topology on $\mathbb{B}\ \subset\ \K({\mathbb{D}})$.
For $\epsilon > 0$, consider $$V = \big\{ (\theta',\phi') : d_{\mathbb{G}} \Big(\big(\gamma', (\theta',\phi')\big), \big( \gamma, (\theta,\phi) \big) \Big) < \epsilon\big\}.$$
Therefore, $$d_{\partial\mathbb{D}}(\theta,\theta') , d_{\partial\mathbb{D}}(\phi,\phi') < \epsilon.$$
Without any loss of generality, we assume that $\epsilon$ is small enough such that the $\epsilon$-tube of the base arcs about $(\theta,\phi)$ doesn't contain any vertex of $\Pi$. The adjoining figure elaborates the schematics.  Consider the open balls $U_1,U_2,...,U_n$ in $\mathbb{D}$ such that $$(\theta,\phi) \subset \cup_{i=1}^{n} U_i,\ (\theta, \phi) \cap U_i \neq \emptyset\ \forall\ i= 1,...,n $$ and each $U_i$ lying inside the $\epsilon$-tube. Since $<U_1,...,U_n>$ is open in $\K({\mathbb{D}})$, therefore $\mathbb{B}\ \cap <U_1,...,U_n>$ is open in $\mathbb{B}$ and is lying in the $\epsilon$-tube. Therefore, we have $$(\theta,\phi) \in \mathbb{B}\ \cap<U_1,...,U_n>\ \subset V.$$

\begin{figure}[h!]
\centering
\includegraphics[width=8.2cm,height=8cm]{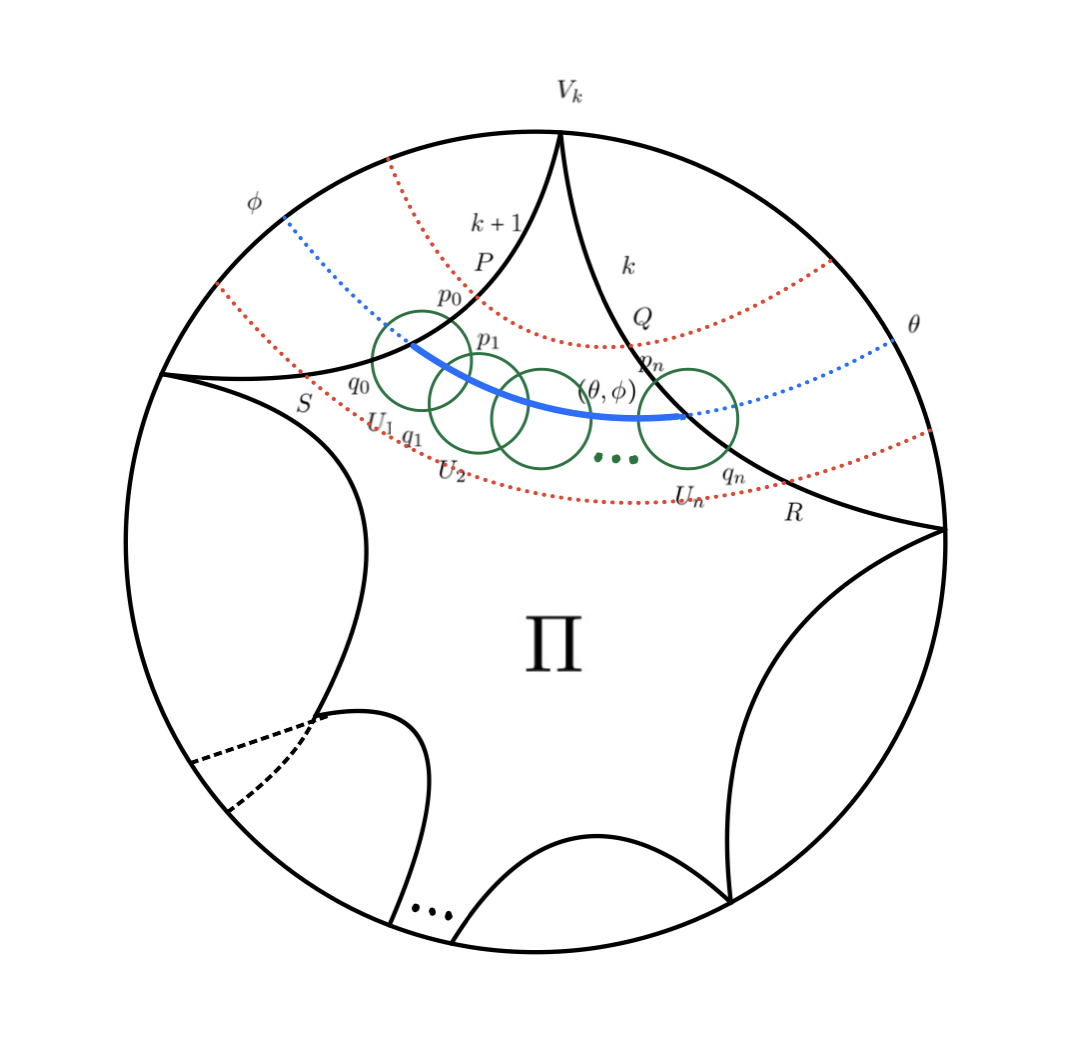}

\caption{An $\epsilon$-tube about a base arc}

\end{figure}

Conversely, consider a basic open set $\mathbb{B}\ \cap <U_1,...,U_n>$ containing a base arc $(\theta, \phi)$.
Without the loss of generality, we assume that $U_i's$ are open discs in $\mathbb{D}$. Define $W_{ij} = \{ p \in \mathbb{D} : p \in U_i\ \cap\ U_j\  \forall i,j \in \{1,...,n\}, i \neq j \}.$ Note that each $W_{ij}$ is either $\emptyset$ or contains two points. Define $W_0 = \{p \in \mathbb{D} : p \in (U_i \cap(\partial \Pi)_k) \cup(U_i \cap(\partial \Pi)_{k+1}) \ \forall\ i = 1,...,n \} .$ Here, $(\partial \Pi)_k$ denotes the side of the polygon with label $k$. Define $W= (\cup_{i,j=1,i \neq j}^{n} W_{ij}) \cup W_0 $ and choose $\delta < inf_{p \in W} (d_{\partial \mathbb{D}}(p,(\theta,\phi)) .$ Then, the $\delta$-tube $V = \{(\theta', \phi') : d_{\mathbb{G}}((\theta',\phi'),(\theta,\phi)) < \delta \}$ lies inside $\mathbb{B}\ \cap <U_1,...,U_n>$.

\eo

\subsection{Ideal Polygons}
\subsubsection{Symbolic Dynamics for Billiards in Ideal Polygons}

Theorem \ref{castle result} ensures that the set $S(\Pi)$ is not dependent on the  choice of the ideal polygon $\Pi$. The elements of $S(\Pi)$ are restricted only by the rules (1) and (2) defined there. 

\bigskip

Let $\Pi$ be an ideal polygon that is symmetrically placed on $\mathbb{D}$, and let $G$ be the space of pointed geodesics on $\Pi$.

Define a map $\tau : \mathbb{G} \to \mathbb{G}$ with its action on $\mathbb{G}$ described as follows :
$$\tau\Big(\big(\gamma, (\theta,\phi) \big)\Big) = \big(\gamma, T(\theta,\phi) \big) \hspace{4mm}\forall\ j \in \mathbb{Z}.$$ 

We will study the dynamics of $(\mathbb{G},\tau)$ under the metric $d_{\mathbb{G}}$.

\bt \label{ideal2} Let $\ \Pi \ \subset \ \mathbb{D}\ be \ an \ ideal \ polygon  \ with \ counter-clockwise\ enumeration  \ 1,...,k$ and $\mathbb{G}$ be the space of pointed geodesics on $\Pi$. Suppose $X$ be the space of all pointed bi-infinite sequences $...a_{-1}.a_0a_1... \in {\{1,....k\}}^{\mathbb{Z}}$ satisfying the rules:\\ (1) $a_j \neq a_{j+1}\ \forall\ j \in \mathbb{Z}$ and \\(2)$...a_{-1}.a_0a_1...$ does not contain an infinitely repeated sequence or bi-infinite sequence of labels of two adjacent sides.\\ Then $(\mathbb{G},\tau) \simeq (X,\sigma)$. \et
\bo Define $h : (\mathbb{G},\tau) \to (X,\sigma)$ by $h\big(\gamma, (\theta,\phi) \big) = ...a_{T^{-1}(\theta,\phi)}.a_{(\theta,\phi)} a_{T(\theta,\phi)}...$

Now $h\big(\gamma, (\theta,\phi) \big)  = h\big(\gamma', (\theta',\phi') \big) 
\Rightarrow ...a_{T^{-1}(\theta,\phi)}.a_{(\theta,\phi)} a_{T(\theta,\phi)}... = ...a_{T^{-1}(\theta',\phi')}.a_{(\theta',\phi')} a_{T(\theta',\phi')} ... $

$\Rightarrow \big(a_{T^n(\theta,\phi)}\big)_{n \in \mathbb{Z}} = \big(a_{T^n(\theta',\phi')}\big)_{n \in \mathbb{Z}}$.

From \cite{castle} , we see that $(T^n(\theta,\phi)\big)_{n \in \mathbb{Z}} = (T^n(\theta',\phi')\big)_{n \in \mathbb{Z}}$ and $ a_{(\theta,\phi)} = a_{(\theta'.\phi')}
\Rightarrow \big(\gamma, (\theta,\phi) \big) = \big(\gamma', (\theta',\phi')\big)$. This gives the injectivity of h.

The surjectivity of $h$ is established from the fact that each $(a_j)_{j \in \mathbb{Z}} \in X$ defines a unique billiard trajectory $\gamma$.
Therefore with corresponding $...a_{-1}.a_0a_1...$, we get a unique base symbol $a_0$, which further picks a base arc $(\theta,\phi)$ on $\gamma$, thereby giving us a unique pointed geodesic in $\mathbb{G}$ i.e.\begin{align*} h\big(\gamma, (\theta,\phi) \big) = ...a_{-1}.a_0a_1... \end{align*}
Thus $ h \circ \tau\Big(\big(\gamma, (\theta,\phi) \big)\Big) = h\bigg(\tau\Big(\big(\gamma, (\theta,\phi) \big)\Big)\bigg) = h\Big(\big(\gamma, T(\theta,\phi) \big)\Big) = h\Big(\big(\gamma, (\theta_1,\phi_1) \big)\Big)$ 

 $ = ...a_{T^{-1}(\theta_1,\phi_1)}.a_{(\theta_1,\phi_1)} a_{T(\theta_1,\phi_1)}...=...a_{T^{-1}T(\theta,\phi)}.a_{T(\theta,\phi)} a_{TT(\theta,\phi)}...= ...a_{(\theta,\phi)}.a_{T(\theta,\phi)} a_{T^2(\theta,\phi)}...$   
 
 $ = \sigma \Big(h\big(\gamma, (\theta,\phi) \big)\Big) = \sigma \circ h \big(\gamma, (\theta,\phi) \big)$.

Therefore, $h \circ \tau = \sigma \circ h$, implying that $h$ is a homomorphism.

Consider an open set $$V = B_{\epsilon}\big(\gamma,(\theta,\phi)\big)$$ in $\mathbb{G}$. Thus $\big(\gamma',(\theta',\phi')\big) \in V $ if and only if $d_{\partial\mathbb{D}}(\theta,\theta'),d_{\partial\mathbb{D}}(\phi,\phi') < \epsilon.$ Tesselate $\mathbb{D}$ with $\Pi$ and its copies generated by reflecting $\Pi$ about its sides and doing the same for the reflected copies along the unfolded geodesic generated by $\gamma$. Let us label the vertices of $\Pi$ in anticlockwise sense by $A_1, A_2,...., A_k$ and the vertices of the $i^{th}$ copy of $\Pi$ by $A^i_1, A^i_2,...., A^i_k$. Define p to be the largest positive integer such that  $$A^i_1, A^i_2,...., A^i_k \not\in (\theta-\epsilon,\theta+\epsilon) \bigtimes (\phi-\epsilon,\phi+\epsilon) \hspace{5mm} \forall\ i = -p,-p+1,...,0,1,...,p.$$
    Then $h^{-1}([x_{-p}...x_{-1}x_0...x_k]) \subseteq V.$ 
    Therefore, $h^{-1}$ is continuous.

\begin{figure}[h!]
\centering
\includegraphics[width=8.5cm,height=7.5cm]{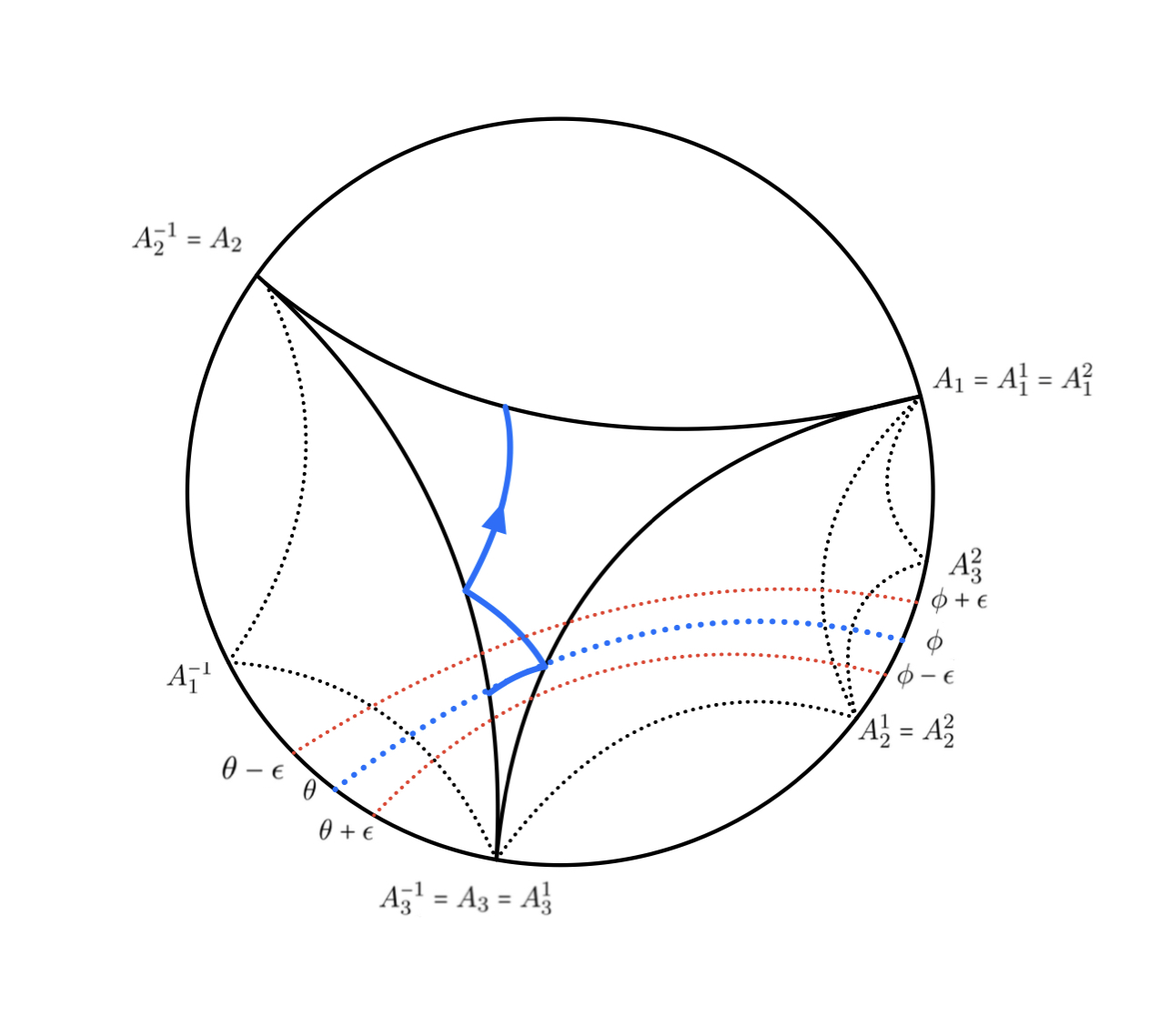}
\caption{$\epsilon$-tube\ about\ an\ unfolded\ billiard\ trajectory}
\end{figure}

Consider an open set $$U = [x_{-m}...x_{-1}.x_0...x_m]$$ in $(X, \sigma)$. Pick an arbitrary  bi-infinite sequence $x \in$ U. Then from corresponding $(x_n)_{n \in \mathbb{Z}}$, we get a billiard trajectory $\gamma$ using \ref{castle result}. By pointing out the base arc $(\theta,\phi)$ corresponding to symbol $x_0$, we get a pointed geodesic $(\gamma,(\theta,\phi))$ whose billiard bi-infinite sequence, we label as $y=(y_i)$. Now, in general $y$ may not be in $U$, but since $x$ and $y$ belong to same equivalence class, there exists an $s$ such that $y_{[s-m,s+m]}$ equals $x_{-m}...x_{-1}.x_0...x_m$. Therefore, $(\gamma,T^{-s}(\theta,\phi))$ has its associated pointed billiard bi-infinite sequence $h(\gamma,T^{-s}(\theta,\phi))\ \in U.$

We construct m future and m past copies of $\Pi$ in $\mathbb{D}$ by reflecting $\Pi$ about its sides as suggested by $$h(\gamma,T^{-s}(\theta,\phi))\ \in U.$$ Label $T^{-s}(\theta,\phi)$ as $({\theta}',{\phi}')$.

 Let $\delta_1$ be defined as follows : \begin{align*} \delta_1 = \displaystyle\min_{i \in \{1,...,k\}}\big\{d_{\partial\mathbb{D}}(A^m_i,\phi'),d_{\partial\mathbb{D}}(A^{-m}_i,\theta')\big\}. \end{align*}
    Choose $\epsilon$ such that $0 < \epsilon <\delta_1$.
    If  $\big(\gamma',(\theta',\phi')\big) \in B_{\epsilon}\big(\gamma,(\theta,\phi)\big),$ then \begin{align*} [h\big(\gamma',(\theta',\phi')\big)]_{[-m,m]} = \displaystyle   x_{-m}...x_{-1}x_0...x_m. \end{align*}
    Thus $ h\big(\gamma',(\theta',\phi')\big) \in U$, i.e. $h\Big(B_{\epsilon}\big(\gamma,(\theta,\phi)\big)\Big) \subseteq U.$ 
    Therefore, $h$ is continuous. \eo

Thus, the space of all  bi-infinite sequences on a k-sided ideal polygon is given by

$X = \{ ...x_{-1}.x_0 x_1...  \in \{1,2,...,k\}^{\mathbb{Z}} : x_i \neq x_{i+1}\ \forall\ i\ and\ ...x_{-1}.x_0 x_1... \neq \overline{ab} ,\ w\overline{ab} ,\ \overline{ab} w \ \text{for\ any\ adjacent}\ a,\ b \in \{1,...,k\}\ \text{and\ word}\ w \}.$

Therefore, $X$ is not closed as the limit points of $X$ of type $\overline{ab}$, $w\overline{ab}$, $\overline{ab} w$ do not lie in $X$. Thus, we further look for the closure of $X$ in $\{1,...,k\}^{\mathbb{Z}}$ and define $\tilde{X} = X \cup X'$ where $X'$ is the set of all limit points of $X$.
Hence, $$\tilde{X} = \{...x_{-1}.x_0 x_1... \in \{1,...,k\}^{\mathbb{Z}} : x_i \neq x_{i+1}\  \forall\ i\}$$ and thereby is an SFT with forbidden set $\{11, 22,...,kk\}$. Thus, $X$ is a dense subset of an SFT. We notice that $\tilde{X}$ is a completion of $X$, therefore it is also a compactification of $X$.

Thus, here we obtain a pair of conjugacies $\pi$ between  $(\mathbb{G},\tau)$ and $(\mathbb{B},T)$ given by $\pi : \mathbb{G} \rightarrow \mathbb{B}$ such that $\pi(\gamma,(\theta,\phi)) = (\theta,\phi)$ and $h$ between $(\mathbb{G},\tau)$ and $(X,\sigma)$ with $h$ defined as above.

In particular, the diagram

$$\begin{CD}
X @>\sigma>> X\\
@V{h^{-1}}VV   @VV{h^{-1}}V\\
\mathbb{G}  @>\tau>> \mathbb{G}\\
@V{\pi}VV   @VV{\pi}V\\
\mathbb{B} @>T>> \mathbb{B}
\end{CD}$$

commutes.

\subsubsection{Dynamical Properties of Billiards in Ideal Polygons}

For an ideal polygon with $k$ vertices, the closure of the space of pointed bi-infinite sequences has forbidden set $$F = \{11,22,...,kk\}.$$ Thus, the adjacency matrix for the corresponding vertex shift is $$A = \begin{pmatrix}0&1&1&. \ .\ .&1\\1&0&1&.\ .\ .&1\\.\\.\\.\\1&1&1&.\ .\ .&0\end{pmatrix}.$$ Its eigenvalues turn out to be $-1,-1,...,-1,k-1$ and further combinatorial computations yield $$B_n(\tilde{X}) = k (k-1)^n.$$ This implies that $$h_{\tilde{X}} = \lim \limits_{n \to \infty} \frac{\log \mathrm{B}_{n}}{n}\ = \lim \limits_{n \to \infty} \frac{\log(k(k-1)^n)}{n}\ = \log(k-1).$$ Since $X$ is a dense subset of $\tilde{X}$, they share the same entropy \cite{hofer}.

\bigskip

We also note that the matrix $A$ is aperiodic and so $(\tilde{X},\sigma)$ is mixing. Since every element of $X'$ consists as subword the infinite strings of the form $ab$ for $a,b \in \{1, \ldots, k\}$, we note that $(X, \sigma)$ is also mixing.

\bigskip

Thus $(\mathbb{G}, \tau)$ is a mixing system with topological entropy $\log(k-1)$, with $k$ being the number of ideal vertices of the polygon $\Pi$.

\subsection{Compact Rational Polygons}
Now, we consider the case when the vertices of $\Pi$ lie inside $\mathbb{D}$. The coding for most bounded systems is non-exact and so is the case here for the compact polygons except for the ones that tile $\mathbb{D}$. Therefore, to ensure that the fundamental domain tesselates the disc, we only consider the polygons with \emph{rational} angles. Recall that, an angle at a vertex of a polygon in $\mathbb{D}$ is called \emph{rational} if it is of the form $\pi / n$ where $n \in \mathbb{N}$ and $n > 1$. We will call the corresponding vertex a \emph{rational vertex}.

\subsubsection{Symbolic Dynamics for Billiards in Compact Rational Polygons}

First, we will establish the coding rules for such rational polygons motivated by the discussion in \cite{ullmo2}. Theorem \ref{compact1} has been originally discussed in \cite{ullmo2}. Here, we present an alternative approach and proof for the same using the basic techniques of topological dynamical systems. We follow this by establishing the  conjugacy between the space of corresponding pointed geodesics and the symbolic space defined under the coding rules presented in Theorem \ref{compact1}.

\bt \label{compact1}
Let $\Pi \subset \mathbb{D}$ be a compact rational polygon with anti-clockwise enumeration of sides labeled $1,2,...,k$. Label the vertices of $\Pi$ as $v_1,...,v_k$ with $\Omega_1,...,\Omega_k$ being the respective
interior angles such that the adjacent sides of $v_i$ are $i$ and $i+1$. Further, assume that $\lambda_i = \pi / \Omega_i  \in \mathbb{N}$ for each $i \in \{1,...,k\}$.
Then an equivalence class of  bi-infinite sequences, $(a_j)$ with $...a_{-1}.a_0a_1... \in \{1,...,k\}^{\mathbb{Z}}$ is in $S(\Pi)$(the space of all equivalence classes of  bi-infinite sequences of $\Pi$) if and only if \\
(1) $...a_{-1}.a_0a_1...$ does not contain any immediate repetitions of symbols i.e., $a_j \neq a_{j+1}\ \forall\ j \in \mathbb{Z}$.\\
(2) $...a_{-1}.a_0a_1...$ does not contain more than $\lambda_i$ repetitions of two successive symbols $i$ and $i+1$ for every $i \in \{1,...,k\}$.\\
Moreover, every equivalence class of such bi-infinite sequences corresponds to one and only one billiard trajectory.
\et

\bo First, we will establish the necessity of (1) and (2). Consider a   bi-infinite sequence $...a_{-1}.a_0a_1... $ attached to a billiard trajectory. In $\mathbb{D}$, two distinct geodesics can have at most one intersection, therefore (1) holds. Now, suppose (2) does not hold for $...a_{-1}.a_0a_1... $ i.e., there exists a subword $w$ of $...a_{-1}.a_0a_1... $ with $\mu_i(> \lambda_i)$ repetitions of letters $i$ and $i+1$. Without any loss of generality, let us assume that $w$ starts with $i$. Then, on unfolding the corresponding part of the trajectory $\mu_i$ times, we get the condition $\Omega_i \mu_i  > \pi$ at the vertex $v_i$, which gives a contradiction.

\smallskip

Now, suppose $...a_{-1}.a_0a_1... $ be a bi-infinite sequence for which (1) and (2) hold. We will construct a unique billiard trajectory defined by it in $\Pi$. We  split our proof into three parts. First we will prove that if we consider a bi-sequence satisfying (1) and (2) which is periodic, it uniquely defines a billiard trajectory. Then, we will show that the set of periodic bi-infinite sequences satisfying (1) and (2) is dense in the set of all bi-infinite sequences satisfying (1) and (2). Lastly, this will allow us to construct a unique billiard trajectory against an arbitrary bi-infinite sequence satisfying (1) and (2) as the limit of a sequence comprising of periodic bi-infinite sequences satisfying (1) and (2).

The set of all periodic  bi-infinite sequences is dense in a full $k$-shift. If we remove bi-infinite sequences containing words $(i \ i+1)^{\mu_i}$ and $(i+1 \ i)^{\mu_i}$ where $\mu_i > \lambda_i = \pi/ \Omega_i\ \forall i$, and the ones for which $a_i=a_{i+1}$ for any $i$, the set of all periodic  bi-infinite sequences in the remaining space is still dense. Indeed, if we start with $...a_{-1}.a_0a_1... $ satisfying (1) and (2), we can define a sequence of periodic points as follows: take $$x^1 = (a_{-1} .a_0 a_1)^\infty ,\ x^2 = (a_{-2} a_{-1}. a_0 a_1 a_2)^\infty$$ and so on. An issue with a typical term of such a sequence can be that the maximal repeating word may have starting and ending letters same or adjacent which violates (1) and (2) respectively. As a remedy to this, whenever such a violation happens, we drop that term from the sequence. This procedure still leaves a subsequence (because of the fact that $...a_{-1}.a_0a_1... $ itself satisfies (2) which comprises of periodic points and converges to $...a_{-1}.a_0a_1... $. Thus, we have the required denseness property.

Suppose $...a_{-1}.a_0a_1... $ is an arbitrary bi-infinite sequence obeying (1) and (2). Let $(x_m)_{m \in \mathbb{N}}$ be a sequence generated by the above construction, then $x^m \rightarrow ...a_{-1}.a_0a_1... $. Now, each $x^m$ has a unique geodesic $\gamma^m$ associated with it. Since, $\mathbb{D}$ is geodesically complete, the sequence of geodesics, $(\gamma^m)_{m \in \mathbb{N}}$ converges to a limit geodesic $\gamma$, which acts as the unique geodesic associated with $...a_{-1}.a_0a_1... $. On folding $\gamma$ back into the fundamental polygon $\Pi$, we get a unique billiard trajectory associated with $[...a_{-1}.a_0a_1...]$.
\eo

\bt \label{compact2} Let $\Pi \subset \mathbb{D}$ be a compact rational polygon with anti-clockwise enumeration of sides labeled $1,2,...,k$. Label the vertices of $\Pi$ as $v_1$,...,$v_k$ with $\Omega_1$,...,$\Omega_k$ being the respective interior angles such that the adjacent sides of $v_i$ are $i$ and $i+1$. Further, assume that $\lambda_i = \pi / \Omega_i  \in \mathbb{N}$ for each $i \in \{1,...,k\}$. Let $\mathbb{G}$ be the space of pointed geodesics on $\Pi$ and $X$ the space of all  bi-infinite sequences $...a_{-1}.a_0a_1... \in \{1,...,k\}^{\mathbb{Z}}$ satisfying (1) and (2) from Theorem \ref{compact1}. Then $(\mathbb{G}, \tau) \simeq  (X,\sigma)$.  \et
\bo
Define $h : (\mathbb{G},\tau) \to (X,\sigma)$ by $h\big(\gamma, (\theta,\phi) \big) = ...a_{T^{-1}(\theta,\phi)}.a_{(\theta,\phi)} a_{T(\theta,\phi)}...$ \\

Now
$h\big(\gamma, (\theta,\phi) \big)  = h\big(\gamma', (\theta',\phi') \big) 
\Rightarrow ...a_{T^{-1}(\theta,\phi)}.a_{(\theta,\phi)} a_{T(\theta,\phi)}... = ...a_{T^{-1}(\theta',\phi')}.a_{(\theta',\phi')} a_{T(\theta',\phi')} ... $

$\Rightarrow \big(a_{T^n(\theta,\phi)}\big)_{n \in \mathbb{Z}} = \big(a_{T^n(\theta',\phi')}\big)_{n \in \mathbb{Z}}$. From Theorem \ref{compact1} , we see that
\begin{align*} (T^n(\theta,\phi)\big)_{n \in \mathbb{Z}} &= (T^n(\theta',\phi')\big)_{n \in \mathbb{Z}} \ \text{and} \ a_{(\theta,\phi)} = a_{(\theta'.\phi')}\\
\Rightarrow \big(\gamma, (\theta,\phi) \big) &= \big(\gamma', (\theta',\phi')\big). \end{align*} This gives the injectivity of h.\\
The surjectivity of $h$ is established from the fact that each $(a_j)_{j \in \mathbb{Z}} \in S(\Pi)$ defines a unique billiard trajectory $\gamma$ as shown in Theorem \ref{compact1}.
Therefore with corresponding $...a_{-1}.a_0a_1...$, we get a unique base symbol $a_0$, which further picks a base arc $(\theta,\phi)$ on $\gamma$, thereby giving us a unique pointed geodesic in $\mathbb{G}$ i.e. $ h\big(\gamma, (\theta,\phi) \big) = ...a_{-1}.a_0a_1... $.

The homomorphism of $h$ follows under the same arguments as in Theorem \ref{ideal2}.

Consider an open set $V = B_{\epsilon}\big(\gamma,(\theta,\phi)\big)$ in $\mathbb{G}$. Thus $\big(\gamma',(\theta',\phi')\big) \in V$ if and only if $$d_{\partial\mathbb{D}}(\theta,\theta'),d_{\partial\mathbb{D}}(\phi,\phi') < \epsilon.$$

 Since the vertices are rational, it allows us to tesselate $\mathbb{D}$ with $\Pi$ and its copies generated by reflecting $\Pi$ about its sides and doing the same for the reflected copies along the unfolded geodesic generated by $\gamma$. Let us label the vertices of $\Pi$ in anticlockwise sense by $A_1, A_2,...., A_k$ and the vertices of $i^{th}$ copy of $\Pi$ by $A^i_1, A^i_2,...., A^i_k$. We note here that as we replicate $\Pi$ along the unfolded geodesic $\gamma$, it shrinks to a point both in future and past in euclidean sense.
Define p to be the largest positive integer such that  $A^i_1, A^i_2,...., A^i_k \not\in (\theta-\epsilon,\theta+\epsilon) \bigtimes (\phi-\epsilon,\phi+\epsilon) \hspace{5mm} \forall\ i = -p,-p+1,...,0,1,...,p.$

 Then $h^{-1}([x_{-p}...x_{-1}x_0...x_k]) \subseteq V.$ Therefore, $h^{-1}$ is continuous.

\begin{figure}[h!]
\centering
\includegraphics[width=7.5cm,height=7.9cm]{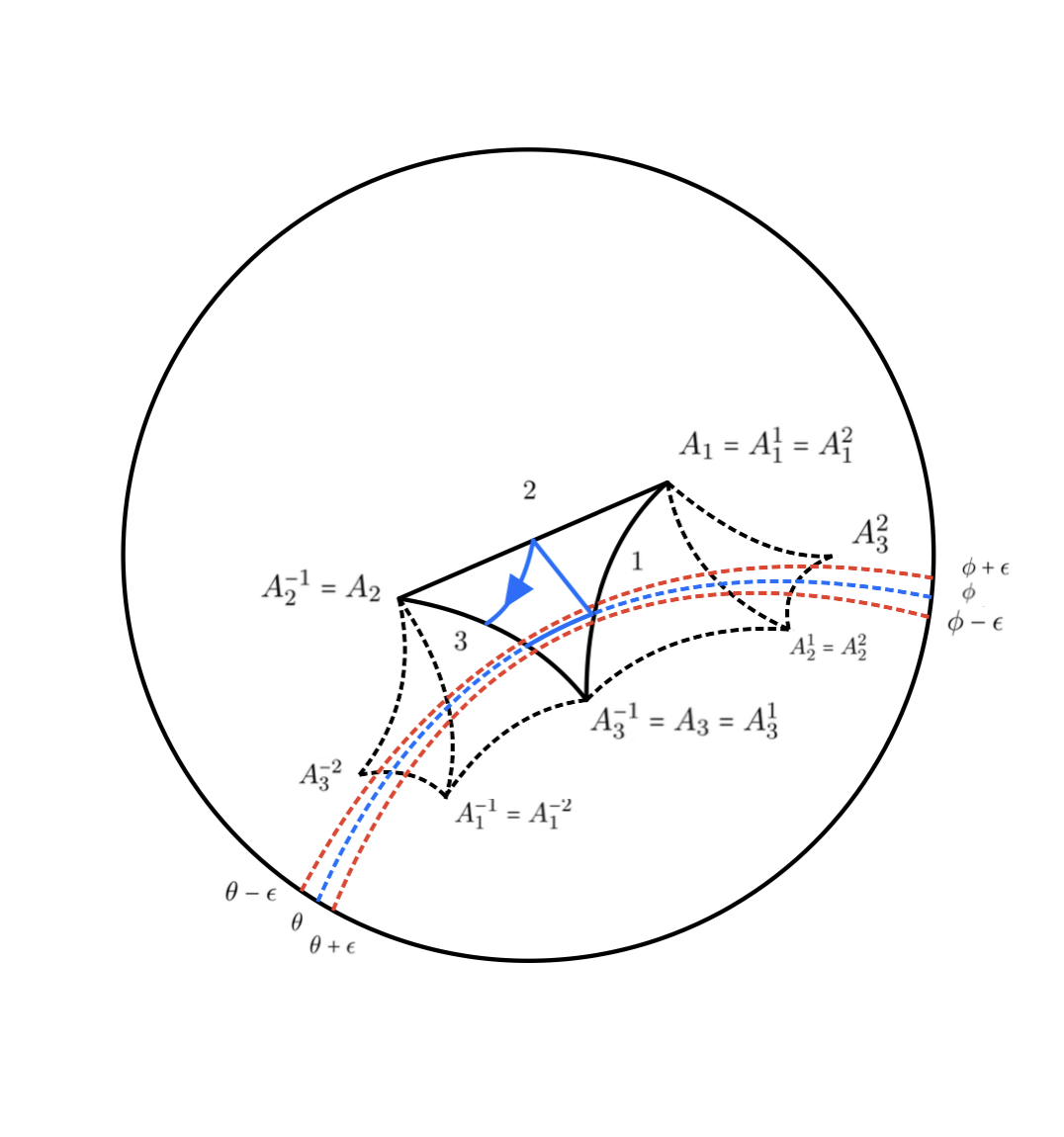}
\caption{$\epsilon$-tube\ about\ an\ unfolded\ billiard\ trajectory}
\end{figure}

Consider an open set $U = [x_{-m}...x_{-1}.x_0...x_m]$ in $(X, \sigma)$. Pick an arbitrary  bi-infinite sequence $x \in$ U. Then from corresponding $(x_n)_{n \in \mathbb{Z}}$, we get a billiard trajectory $\gamma$ using Theorem \ref{compact1}. By pointing out the base arc $(\theta,\phi)$ corresponding to symbol $x_0$, we get a pointed geodesic $(\gamma,(\theta,\phi))$ that corresponds to a bi-infinite sequence which we label as $y=(y_i)$. Now, in general $y$ may not be in $U$, but since $x$ and $y$ belong to same equivalence class, there exists an $s$ such that $y_{[s-m,s+m]}$ equals $x_{-m}...x_{-1}.x_0...x_m$. Therefore, $(\gamma,T^{-s}(\theta,\phi))$ has its associated  billiard bi-infinite sequence $h(\gamma,T^{-s}(\theta,\phi))\ \in U.$

We construct $ m $ future and $ m $ past copies of $\Pi$ in $\mathbb{D}$ by reflecting $\Pi$ about its sides as suggested by $h(\gamma,T^{-s}(\theta,\phi))\ \in U$. Label $T^{-s}(\theta,\phi)$ as $({\theta}',{\phi}')$.

 Let $\delta_1$ be defined as follows : \begin{align*} \delta_1 = \displaystyle\min_{i \in \{1,...,k\}}\big\{d_{\partial\mathbb{D}}(A^m_i,\phi'),d_{\partial\mathbb{D}}(A^{-m}_i,\theta')\big\}. \end{align*}
    Choose $\epsilon$ such that $0 < \epsilon <\delta_1$.
    If  $\big(\gamma',(\theta',\phi')\big) \in B_{\epsilon}\big(\gamma,(\theta,\phi)\big)$,
    then \begin{align*} [h\big(\gamma',(\theta',\phi')\big)]_{[-m,m]} = \displaystyle x_{-m}...x_{-1}x_0...x_m. \end{align*}
    Thus $ h\big(\gamma',(\theta',\phi')\big) \in U$ , i.e. $h\Big(B_{\epsilon}\big(\gamma,(\theta,\phi)\big)\Big) \subseteq U.$
    Therefore, $h$ is continuous.
    
    \smallskip
    
    Thus, $(\mathbb{G}, \tau) \simeq  (X,\sigma)$.
\eo

\br We note that both $S(\Pi)$ and $X$ generate the same language. $S(\Pi)$ gives the trajectories whereas $X$ gives each point in the trajectory.\er

\subsubsection{Dynamical Properties of Billiards in Compact Rational Polygons}

For a compact polygon with $k$ vertices, the space of pointed bi-infinite sequences has forbidden set $$\cF = \{11,22,...,kk,{\underbrace{1212\cdots }_{{(1+\lambda_1)-times}}},  {\underbrace{2121\cdots }_{{(1+\lambda_1)-times}}} , {\underbrace{2323\cdots }_{{(1+\lambda_2)-times}}} , {\underbrace{3232\cdots }_{{(1+\lambda_2)-times}}} , ..., {\underbrace{1k1k\cdots }_{{(1+\lambda_k)-times}}}, {\underbrace{k1k1\cdots }_{{(1+\lambda_k)-times}}} \} .$$ 

\bex We consider a particular case here where we consider a compact triangle with $\lambda_1 =\lambda_2 = \lambda_3 = 4$. Here, we can declare $X$ to be a $4-step\ SFT$ with

 $$\cF_{4+1} = \{ {\underbrace{wiiw' }_{{5-times}}}\ \forall\ i,\ \forall w,\ \forall w',\ 12121,\ 21212,\ 23232,\  32323,\ 13131,\  31313\}.$$
Since $X$ is a $4-step\ SFT$, we go for $X^{[4+1]} = X_G$.

Note that the associated transition matrix $A$ here is aperiodic, and so $ ({X},\sigma)$ is mixing. Also this $A$ has a positive Perron eigenvalue which gives a positive topological entropy for $  ({X},\sigma)$. 

\smallskip

Thus here $(\mathbb{G}, \tau)$ is a mixing system with positive topological entropy. \eex

In general both $(\mathbb{G}, \tau)$ and  $(X,\sigma)$ will be mixing systems, that has been proved later in Theorem \ref{mixing}. And so will have positive topological entropy.

\subsection{Semi-Ideal Rational  Polygons}

After establishing the codes for the billiards on ideal and compact rational polygons, we now do the same for the case where some vertices of the polygon $\Pi$ sit on $\partial \mathbb{D}$ and some in $\mathbb{D}$ such that these polygons also tile $\mathbb{D}$. We call such polygons \emph{semi-ideal rational}.   The coding rules for the billiards on such polygons is the natural amalgamation of the rules from the ideal and the compact rational polygon case.

\bt \label{semiideal1} Let $\Pi \subset \mathbb{D}$ be a semi-ideal rational polygon with anti-clockwise enumeration of sides labeled $1,2,...,k$. Label the vertices of $\Pi$ as $v_1$,...,$v_k$ with $\Omega_1$,...,$\Omega_k$ being the respective interior angles such that the adjacent sides of $v_i$ are $i$ and $i+1$. Further, assume that $v_i \in \mathbb{D}\ \forall\ i \in \Lambda \subset \{1,...,k\}$ and $v_i \in \partial \mathbb{D}\ \forall\ i \in \{1,...,k\} - \Lambda$ with ${\lambda}_i = \pi / {\Omega}_i \in\ \mathbb{N}\ \forall\ i \in \Lambda$. Then an equivalence class $[...a_{-1}.a_0a_1... ]$ with $ ...a_{-1}.a_0a_1... \in \{1,...,k\}^{\mathbb{Z}}$ is in $S(\Pi)$(the space of all equivalence classes of pointed billiard bi-infinite sequences of $\Pi$) if and only if \\
(1) $...a_{-1}.a_0a_1... $ does not contain any immediate repetitions of symbols i.e., $a_j \neq a_{j+1}\ \forall\ j \in \mathbb{Z}$.\\
(2) $...a_{-1}.a_0a_1... $ does not contain more than $\lambda_i$ repetitions of two successive symbols $i$ and $i+1$ for every $i \in \Lambda$.\\
(3) $...a_{-1}.a_0a_1... $ does not contain an infinitely repeated sequence or bi-sequence of labels of two adjacent sides $i$ and $i+1$ $\forall\ i \in \{1,...,k\}- \Lambda$.\\
Moreover, every equivalence class of such billiard induced bi-infinite sequences corresponds to one and only one billiard trajectory.\et

\bo
Working on the similar lines, as for the case of ideal and compact case, we will start with establishing the necessity of (1), (2) and (3). Consider a  bi-infinite sequence $...a_{-1}.a_0a_1... $ attached to a billiard trajectory. Since, two geodesics in a hyperbolic plane can have at most one intersection, therefore (1) holds. Now, suppose (2) does not hold for $...a_{-1}.a_0a_1...$ with $$[...a_{-1}.a_0a_1...] \in S(\Pi).$$ This means that for some $i \in \Lambda$, there exists a subword $w$ of $...a_{-1}.a_0a_1... $ with $\mu_i( > \lambda_i)$ repetitions of letters $i$ and $i+1$. Following the same line of argument as in Theorem \ref{compact1}, we get a contradiction. Suppose (3) does not hold for $...a_{-1}.a_0a_1...$ with $[...a_{-1}.a_0a_1...] \in S(\Pi).$

 This means that for some $i \in \{1,...k\}-\Lambda$, it contains an infinitely repeated sequence or bi-sequence of labels of sides $i$ and $i+1$. Under Cayley transformation, we can shift the whole billiard table to $\mathbb{H}$ and an appropriate isometry of $\mathbb{H}$ allows us to assume that the two sides $i$ and $i+1$ are represented on $\mathbb{H}$ by the vertical lines $x=0$ and $x=1$. Now, our trajectory hits these two sides repeatedly infinite number of times. if we unfold this part of the trajectory in $\mathbb{H}$, the unfolded part of the trajectory being a geodesic lies on a semi-circle centered at real axis. Since a semi-circle can intersect only  finitely many lines from the family $\{x=n : n \in \mathbb{N}\}$, therefore we get a contradiction to our assumption that the trajectory has infinitely many parts consisting of consecutive hits on sides $i$ and $i+1$. Thus, the necessity of conditions (1), (2) and (3) is established. Now suppose $...a_{-1}.a_0a_1...$ be a  bi-infinite sequence for which (1), (2) and (3) hold. We will construct a unique billiard trajectory defined by it in $\Pi$. We structure the rest of the proof on the similar lines as in Theorem \ref{compact1}, going through the three main stages. First, we will prove that if we consider a bi-infinite sequence satisfying (1), (2) and (3) which is periodic, it uniquely defines a billiard trajectory. Then, we will show that the set of periodic bi-infinite sequences satisfying (1), (2) and (3) is dense in the set of all bi-infinite sequences satisfying (1), (2) and (3). Lastly, this will allow us to construct a unique billiard trajectory against an arbitrary bi-infinite sequence satisfying (1), (2) and (3) as the limit of a sequence comprising of periodic bi-infinite sequences satisfying (1), (2) and (3).

The set of all periodic  bi-infinite sequences is dense in a full k-shift. If we remove  bi-infinite sequences containing words $(i \ i+1)^{\mu_i}$ and $(i+1 \ i)^{\mu_i}$ where $\mu_i > \lambda_i = \pi/ \Omega_i$  for some $i \in \Lambda$, or containing an infinitely repeated sequence or bi-infinite sequence of labels of two adjacent sides $i$ and $i+1$ for some $ i \in \{1,...,k\}- \Lambda$, or the ones for which $a_i=a_{i+1}$ for any $i$, the set of all periodic bi-infinite sequences in the remaining space is still dense. Indeed, if we start with $...a_{-1}.a_0a_1... $ satisfying (1), (2) and (3), we can define a sequence of periodic points as follows: take $$x^1 = (a_{-1} .a_0 a_1)^\infty ,\ x^2 = (a_{-2} a_{-1}. a_0 a_1 a_2)^\infty$$ and so on. An issue with a typical term of such a sequence can be that the maximal repeating word may have starting and ending letters same or adjacent which violates (1), (2) and (3) respectively. As a remedy to this, whenever such a violation happens, we drop that term from the sequence. This procedure still leaves a subsequence (because of the fact that $...a_{-1}.a_0a_1... $ itself satisfies (2) and (3) which comprises of periodic points and converges to $...a_{-1}.a_0a_1... $. Thus, we have the required denseness property.

Suppose $...a_{-1}.a_0a_1... $ is an arbitrary bi-infinite sequence obeying (1), (2) and (3). Let $(x_m)_{m \in \mathbb{N}}$ be a sequence generated by the above construction, then $x^m \rightarrow ...a_{-1}.a_0a_1... $. Now, each $x^m$ has a unique geodesic $\gamma^m$ associated with it. Since, $\mathbb{D}$ is geodesically complete, the sequence of geodesics, $(\gamma^m)_{m \in \mathbb{N}}$ converges to a limit geodesic $\gamma$, which acts as the unique geodesic associated with $...a_{-1}.a_0a_1... $. On folding $\gamma$ back into the fundamental polygon $\Pi$, we get a unique billiard trajectory associated with $[...a_{-1}.a_0a_1...]$.
\eo

\bt \label{semi-ideal2} Let $\Pi \subset \mathbb{D}$ be a semi-ideal rational polygon with anti-clockwise enumeration of sides labeled $1,2,...,k$. Label the vertices of $\Pi$ as $v_1$,...,$v_k$ with $\Omega_1$,...,$\Omega_k$ being the respective interior angles such that the adjacent sides of $v_i$ are $i$ and $i+1$. Further, assume that $v_i \in \mathbb{D}\ \forall\ i \in \Lambda \subset \{1,...,k\}$ and $v_i \in \partial \mathbb{D}\ \forall\ i \in \{1,...,k\} - \Lambda$ with ${\lambda}_i = \pi / {\Omega}_i \in\ \mathbb{N}\ \forall\ i \in \Lambda$. Let $\mathbb{G}$ be the space of pointed geodesics on $\Pi$ and $X$ the corresponding space of all  bi-infinite sequences $...a_{-1}.a_0a_1... \in \{1,...,k\}^{\mathbb{Z}}$ satisfying (1), (2) and (3) from Theorem \ref{semiideal1}. Then $(\mathbb{G}, \tau) \simeq  (X,\sigma)$.  \et

\bo
Define $h : (\mathbb{G},\tau) \to (X,\sigma)$ by $h\big(\gamma, (\theta,\phi) \big) = ...a_{T^{-1}(\theta,\phi)}.a_{(\theta,\phi)} a_{T(\theta,\phi)}...$ 

Now
$h\big(\gamma, (\theta,\phi) \big)  = h\big(\gamma', (\theta',\phi') \big) 
\Rightarrow ...a_{T^{-1}(\theta,\phi)}.a_{(\theta,\phi)} a_{T(\theta,\phi)}... = ...a_{T^{-1}(\theta',\phi')}.a_{(\theta',\phi')} a_{T(\theta',\phi')} ... $

$\Rightarrow \big(a_{T^n(\theta,\phi)}\big)_{n \in \mathbb{Z}} = \big(a_{T^n(\theta',\phi')}\big)_{n \in \mathbb{Z}}$. From Theorem \ref{semiideal1} , we see that $(T^n(\theta,\phi)\big)_{n \in \mathbb{Z}} = (T^n(\theta',\phi')\big)_{n \in \mathbb{Z}}$ and $a_{(\theta,\phi)} = a_{(\theta'.\phi')} \Rightarrow \big(\gamma, (\theta,\phi) \big) = \big(\gamma', (\theta',\phi')\big)$.  This gives the injectivity of h.

The surjectivity of $h$ is established from the fact that each $(a_j)_{j \in \mathbb{Z}} \in S(\Pi)$ defines a unique billiard trajectory $\gamma$ as shown in Theorem \ref{semiideal1}.
Therefore with corresponding $...a_{-1}.a_0a_1...$, we get a unique base symbol $a_0$, which further picks a base arc $(\theta,\phi)$ on $\gamma$, thereby giving us a unique pointed geodesic in $\mathbb{G}$ i.e.  $h\big(\gamma, (\theta,\phi) \big) = ...a_{-1}.a_0a_1... $.

The homomorphism of $h$ follows under the same arguments as in Theorem \ref{ideal2}.

Consider an open set $V = B_{\epsilon}\big(\gamma,(\theta,\phi)\big)$ in $\mathbb{G}$. Thus $\big(\gamma',(\theta',\phi')\big) \in V$ if and only if $$d_{\partial\mathbb{D}}(\theta,\theta'),d_{\partial\mathbb{D}}(\phi,\phi') < \epsilon.$$ Since the vertices are rational or ideal, it allows us to tesselate $\mathbb{D}$ with $\Pi$ and its copies generated by reflecting $\Pi$ about its sides and doing the same for the reflected copies along the unfolded geodesic generated by $\gamma$. Let us label the vertices of $\Pi$ in anticlockwise sense by $A_1, A_2,...., A_k$ and the vertices of $i^{th}$ copy of $\Pi$ by $A^i_1, A^i_2,...., A^i_k$. We note here that as we replicate $\Pi$ along the unfolded geodesic $\gamma$, it shrinks to a point both in future and past in euclidean sense. The vertices that lie on $\partial\mathbb{D}$ stay put, whereas the ones lying in $\mathbb{D}$ move towards the boundary under both future and past limits.
Define p to be the largest positive integer such that  $$A^i_1, A^i_2,...., A^i_k \not\in (\theta-\epsilon,\theta+\epsilon) \bigtimes (\phi-\epsilon,\phi+\epsilon) \hspace{5mm} \forall\ i = -p,-p+1,...,0,1,...,p.$$ Then $h^{-1}([x_{-p}...x_{-1}x_0...x_k]) \subseteq V.$ Therefore, $h^{-1}$ is continuous.

\begin{figure}[h!]
\centering
\includegraphics[width=8.6cm,height=8cm]{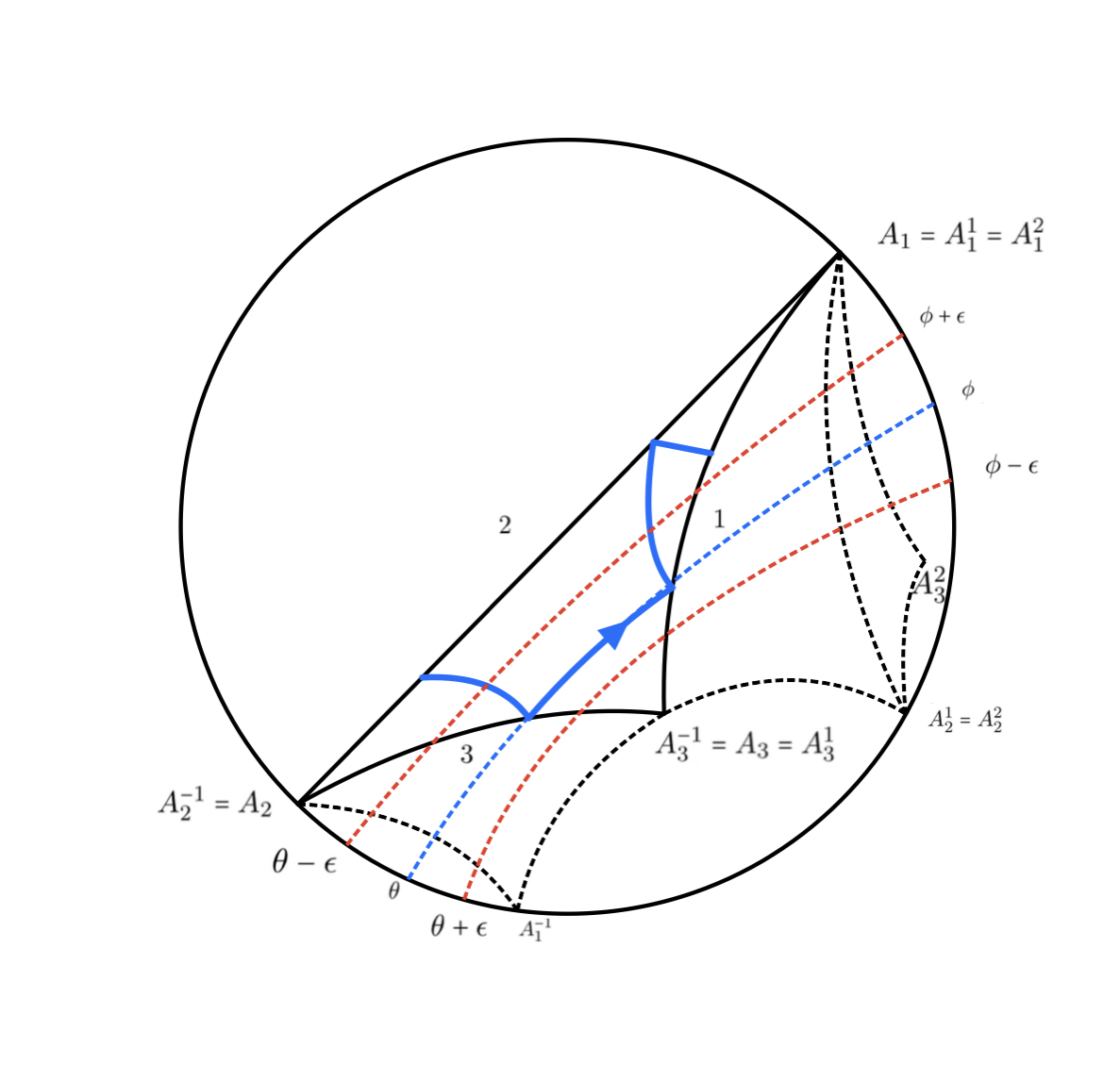}
\caption{$\epsilon$-tube\ about\ an\ unfolded\ billiard\ trajectory}
\end{figure}

Consider an open set $$U = [x_{-m}...x_{-1}.x_0...x_m]$$ in $(X, \sigma)$. Pick an arbitrary  bi-infinite sequence $x \in$ U. Then from corresponding $(x_n)_{n \in \mathbb{Z}}$, we get a billiard trajectory $\gamma$ using Theorem \ref{semiideal1}. By pointing out the base arc $(\theta,\phi)$ corresponding to symbol $x_0$, we get a pointed geodesic $(\gamma,(\theta,\phi))$ corresponding to the bi-infinite sequence which  we label as $y=(y_i)$. Now, in general $y$ may not be in $U$, but since $x$ and $y$ belong to same equivalence class, $\exists$ an $s$ such that $y_{[s-m,s+m]}$ equals $x_{-m}...x_{-1}.x_0...x_m$. Therefore, $(\gamma,T^{-s}(\theta,\phi))$ has its associated  billiard bi-infinite sequence $h(\gamma,T^{-s}(\theta,\phi))\ \in U$.
We construct m future and m past copies of $\Pi$ in $\mathbb{D}$ by reflecting $\Pi$ about its sides as suggested by $$h(\gamma,T^{-s}(\theta,\phi))\ \in U.$$ Label $T^{-s}(\theta,\phi)$ as $({\theta}',{\phi}')$. Let $\delta_1$ be defined as follows : \begin{align*} \delta_1 = \displaystyle\min_{i \in \{1,...,k\}}\big\{d_{\partial\mathbb{D}}(A^m_i,\phi'),d_{\partial\mathbb{D}}(A^{-m}_i,\theta')\big\}. \end{align*}
    Choose $\epsilon$ such that $0 < \epsilon <\delta_1$.
    If  $\big(\gamma',(\theta',\phi')\big) \in B_{\epsilon}\big(\gamma,(\theta,\phi)\big),$
    then $ [h\big(\gamma',(\theta',\phi')\big)]_{[-m,m]} = \displaystyle x_{-m}...x_{-1}x_0...x_m$. 
    
    Thus $ h\big(\gamma',(\theta',\phi')\big) \in U $, i.e. $h\Big(B_{\epsilon}\big(\gamma,(\theta,\phi)\big)\Big) \subseteq U.$
    Therefore, $h$ is continuous.
    
    \smallskip
    
    Thus $(\mathbb{G}, \tau) \simeq  (X,\sigma)$.
\eo

We note that if $\tilde{X}$ is the closure of $X$ in $\{ 1,2, \ldots, k\}^{\Z}$, then $(\tilde{X}, \sigma)$ is an SFT.

\subsubsection{Dynamical Properties of Billiards in Semi-Ideal Rational Polygons}

We note that in case of rational polygons the SFT $(\tilde{X}, \sigma)$ is mixing. We present a proof for the semi-ideal case and note that it works the same for the compact case(where $\tilde{X} = X$).

\bt \label{mixing} Let $\Pi \subset \mathbb{D}$ be a semi-ideal rational polygon with anti-clockwise enumeration of sides labeled $1,2,...,k$. Label the vertices of $\Pi$ as $v_1$,...,$v_k$ with $\Omega_1$,...,$\Omega_k$ being the respective interior angles such that the adjacent sides of $v_i$ are $i$ and $i+1$. Further, assume that $v_i \in \mathbb{D}\ \forall\ i \in \Lambda \subset \{1,...,k\}$ and $v_i \in \partial \mathbb{D}\ \forall\ i \in \{1,...,k\} - \Lambda$ with ${\lambda}_i = \pi / {\Omega}_i \in\ \mathbb{N}\ \forall\ i \in \Lambda$. Let $\mathbb{G}$ be the space of pointed geodesics on $\Pi$ and $X$ the corresponding space of codes $...a_{-1}.a_0a_1... \in \{1,...,k\}^{\mathbb{Z}}$ satisfying (1), (2) and (3) from \ref{semiideal1}. Then $(\tilde{X},\sigma)$ is mixing.  \et

\bo
We start with assumption that $k \geq 6$. For every pair of words $u,v \in \cL(\tilde{X})$, our aim is to declare an $N \in \mathbb{N}$ such that a   $w \in \cL(\tilde{X}) $ can be chosen with $|w| = n$ for all $n \geq N$ with $uwv \in \cL(\tilde{X})$. Suppose the word $u$ ends in $ri$ and the word $v$ starts with $jl$, then we can choose two distict labels $p, q$ distinct from $r,i,j,l$. We construct $w$ as $(pq)^d i(pq)^d i ...(pq)^d$, where the $(pq)$ blocks are repeated $c$ times. Here, $c$ can be chosen arbitrarily large but $d$ has a natural restriction for adjacent $p,q$ under the coding rules that we have laid down in Theorem \ref{semi-ideal2}. This establishes the mixing property for any $\tilde{X}$ associated with a semi-ideal rational polygon with sides equal or more than $6$.

Next, we consider the case where $k = 5$. Under the same notation as above, suppose all of $r,i,j,l$ are distinct. Then, we construct $w$ as $(pj)^di (pj)^di...(pj)^di$ by choosing $p $ distinct from $r,i,j,l$. If two of $r,i,j,l$ are equal, say $i=l \neq j,r$, then we choose $p,q \neq i,l,r$ and construct $w$ as $(pq)^di(pq)^di...(pq)^d$. Similar constructions can be done for $w$ is any other similarities occur among $r,i,j,l$ as we get more room to wiggle around. This establishes the mixing property for any $\tilde{X}$ associated with a semi-ideal rational polygon with $5$ sides.

Now, we consider the case $k= 4$. The situations where more than one similarities occur among $r,i,j,l$ can be handled in same way as discussed above. If all of $r,i,j,l$ are distinct, then we construct $w$ as $(jl)^di(jl)^di...(jl)^di(ri)^e$ where $d,e$ are chosen respecting the respective coding rules. If exactly two of $r,i,j,l$ are equal, say $i=l$, then we choose $p \neq i,j,r$ and construct $w = (pr)^di(pr)^di...(pr)^di(pr)^d$. Other cases of similarity among $r,i,j,l$ can be handled in the same way. This establishes the mixing property for any $\tilde{X}$ associated with a semi-ideal rational polygon with $4$ sides.

Lastly, for $k=3$ atleast two of the labels $r,i,j,l$ are same. If exactly two are same, say $i=j$, then we construct $w$ as $(lr)^di(lr)^di...(lr)^d$. Similarly, if $i=l$, then we construct $w$ as $(jr)^di(jr)^di...i(jr)^d$. If two matchings occur among $r,i,j,l$, say $j=r,\ l=i$, then we choose $p \neq i,r$ and construct  $w$ as $(pr)^di(pr)^d...i(pr)^di(rp)^d$. If more than two matchings occur among $r,i,j,l$, it allows us to choose distinct $p,q$ from $r,i,j,l$ as we get more room to wiggle around and similar constructions for $w$ follow as discussed above. This establishes the mixing property for any $\tilde{X}$ associated with a semi-ideal rational polygon with $3$ sides.

\eo

\br Note that this also gives mixing for $(X, \sigma)$ in the compact rational case. In case of semi-ideal rational polygons, the set $\tilde{X} \setminus X$ consists of  strings of the form $i (i+1)$ as described in Theorem \ref{semiideal1}. This implies that $(X, \sigma)$ is also mixing.

Thus $(\mathbb{G}, \tau)$ is mixing. \er

This also gives positive topological entropy for $(\mathbb{G}, \tau)$.

\section{Convergence in the Space of all Subshifts and Hyperbolic Polygons}

\subsection{The Space of Subshifts}

Let $\mathfrak{S}$ be the set of all shift spaces of bi-infinite sequences defined over a finite alphabet with cardinality at least $2$.

\bigskip

For $k \geq 2$, let $\A_k = \{1,2, \ldots, k\}$. Define $X_k = \A_k^{\Z}$, for each $k \geq 2$, the space of all  bi-infinite sequences taking values in $\{1,2,\ldots ,k\}$ indexed by $\Z$, endowed with the product topology. Then, for each $k \geq 2$, $X_k$ is metrizable with the metric defined in Equation (\ref{metric}), and can be taken as a compact metric  space. Thus for each $k \geq 2$, we have subshifts  $(X_k, \sigma)$.

\bigskip

Now each $X_k \subsetneq \N^{\Z}$, $k \geq 2$. Any  shift space $Y \in \fS$  can be considered isometric to a subshift $X \subseteq
X_k$ for some $k \geq 2$. Thus with respect to the Gromov-Hausdorff metric as defined in Equation (\ref{2}), every subshift in $\fS$ is isometrically imbedded in $\N^{\Z}$. Thus $\N^{\Z}$ can be considered as the universal Urysohn space for all subshifts of bi-infinite sequences over a finite alphabet. Every such shift space $Y \in \fS$ is isometric to a compact subset of $\N^{\Z}$. Hence we can consider $\fS \subseteq \K({\N^{\Z}})$.

\bigskip

Thus each subshift $(Y,\sigma)$, for $Y \in \fS$, can be considered as a subsystem of $(\N^{\Z}, \sigma)$. We hence call $(\N^{\Z}, \sigma)$ as \emph{the universal shift space}, and every subshift from $\fS$  is its subsystem.

\bigskip

Consider the metric space $(\K({\N^{\Z}}), d_H)$. Thus $d_H$ gives a metric on $\fS$ and we lift that metric on the set of subshifts $(Y, \s)$ for $Y \in \fS$. We note that this lift is in the sense of Gromov since we identify isometric shift spaces in $\fS$ with some subset of $X_k$ for some $k$.

\bigskip

Our special emphasis will be on $\mathfrak{sft} \subseteq \fS$, where $\sft$ is the
subclass of all subshifts of finite type in $\fS$. We see that  $\sft$ is a very important subset of $\fS$.

\bl For every distinct  $X,Y \in \sft$, $d_H(X,Y) > 0$. \el

\bo Suppose that there exists no $\delta >0$ such that $d_H(X,Y) > \delta$. Then for every $x \in X$, there exists a $y \in Y$ such that $d(x,y) < \delta$. This means that $ \ \exists \ k>0$ such that $x_{[-k,k]} = y_{[-k,k]}$. Taking $ \delta \to 0$, this means that both $x,y$ have the same words of length $2k+1$, for large $k$. This implies $x = y$, i.e. we get $X = Y$.\eo

\bt $\sft$ is a dense subset of $\fS$. \et

\bo   Let $X \in \fS \setminus \sft$ be a shift space and let   $\cF$ be an infinite collection of blocks such that $X = X_{\cF}$.

We borrow this idea from \cite{marcus}. The elements of $\cF$ can be arranged lexicographically, and accordingly let $\cF = \{f_1, f_2, \ldots, f_n, \ldots \}$.

Define $\cF_j = \{f_1, \ldots, f_j\}$ for $j \geq 1$. Then
$$X = \bigcap \limits_{j \in \N} X_{\cF_j}.$$

Thus, $X_{\cF_j} \to X$. \eo

\smallskip

We try to understand the space $\fS \subseteq \K({\N^{\Z}})$.

Suppose $X_n \to X$ in $\fS$. Then there exists a sequence $\{k_n\} \nearrow \infty$ such that

$$ \forall x \in X \ \exists \ x_n \in X_n \ \text{such that} \ x_{[-k_n, k_n ]} = {x_n}_{[-k_n, k_n]}.$$

\smallskip

We note that the set of finite subsets of $\N^{\Z}$ is dense in $ \K({\N^{\Z}})$. But all of them are not contained in $\fS$. We consider $Z \in \fS$ such that $|Z| < \infty$. Since $\s(Z) \subset Z$, and $\s$ is a homeomorphism, $\s$ just permutes the elements of $Z$. Thus, $Z$ is just a periodic orbit or a finite union of periodic orbits. Also, this $Z$ can be considered as a vertex shift over a finite alphabet and so $Z \in \sft$.

Again we note that for every periodic orbit  $\cO(x) \in  \N^{\Z}$, there exists no other element of $\fS$ that can be contained in $[\cO(x)]_\ep$ for every  $\ep > 0$. Thus each finite element of $\sft$ is isolated.

Let $X \in \sft$ be infinite. Then there is a finite $\cF$ such that $X = X_{\cF}$. Let $\A \subset \N$ be the alphabet set for $X$. Then there exists $a,b \in \A, \ a \neq b$ with $ab^N \notin \cF$ for some $N \in \N$. Else $X$ will be finite.

Define $\cF_j = \cF \cup \{ ab^{N +j} \}$. Then $X_{\cF_1} \subsetneqq X_{\cF_2} \subsetneqq \ldots \subsetneqq X_{\cF_j} \subsetneqq \dots$ and
$$X = \overline{\bigcup \limits_{j \in \N} X_{\cF_j}}.$$

Hence $X_{\cF_j} \to X$.

This gives the complete topological structure of the space $\fS$. Thus, we can say:

\bt Every finite element of $\fS$ is isolated and every infinite element is an accumulation point in $\fS$.

Further, for every infinite $X \in \fS$ there exists a sequence $\{X_{\cF_j}\}$ in $\sft$ such that
$X_{\cF_j} \to X$.\et

\subsection{Convergence of Dynamics in Subshifts}

Let us look into the statement `` $X_n \to X$ in $\fS$''.

\bd We say that the subshift $(X,\s)$ is the \emph{limit of subshifts} $\{(X_n, \s)\}$ and write $\lim \limits_{n \to \infty} \ (X_n, \s) = (X,\s)$ in $(\N^{\Z}, \s)$ if and only if $X_n \to X$ in $\fS$. \ed

\br We note that $(X_n, \s) \to (X,\s) \ \cong \ X_n \to X$.\er

We compare the dynamical properties of the sequence $\{(X_n, \s)\}$ of subshifts in $\sft$ to those of the limit $(X,\s)$. We look into the dynamical properties of the limit $(X,\s)$ that are inherited by this limiting process.

\br We note that as mentioned in \cite{marcus}, if $X_n \to X$ then the entropies $h(X_n) \to h(X)$. By the observation in the previous section, for every $X \in \fS$, we can choose a sequence $\{X_n\}$ in $\sft$ with $X_n \to X$. Now $h(X_n) = \log \lambda_n$ where $\lambda_n$ is the Perron eigenvalue corresponding to $X_n \in \sft$, for every $n$. Thus, each $h(X)$ can be considered as a limit of a sequence of the form $\log \lambda_n$.

This also implies that for $X \in \fS$, if $X_n \to X$ with $\liminf \limits_{n \to \infty} h(X_n) > 0$ then $h(X)>0$. \er

\bpr If $X = \overline{\cO(x)}$ for some $x \in \N^{\Z}$, then there exists a sequence  $\{X_n\}$  in $\sft$ with $|X_n| < \infty$ such that $X_n \to X$. \epr

\bo For some $K > 0$ define succesively

$u_1 = x_{[-K,K]}, \ u_2 = x_{[-K-1,K+1]}, \ \ldots, \ u_n = x_{[-K-n,K+n]}, \ \ldots$.

Define $X_n = \cO(u_n^\infty)$, for $j \geq 1$, where $a^\infty = \ldots aaaaaaa \ldots$ the infinite concatenation of $a$ with itself. Since each $X_n$ is finite, it is an SFT.

Now $X = \overline{\cO(x)}$ and so  $X_n \to X$. \eo

\br We note that when $(X,\s)$ is minimal, for some $X \in \fS$, $X$ is obtained as a limit of a sequence of finite SFTs. Now if $X$ has no isolated points then $X$ is sensitive. Since a finite SFT is always equicontinuous, we note that equicontinuity is not preserved by taking limits. However,  recurrence is preserved on taking limits. \er

We can say more:

\bt For   $X \in \fS$, if $X_n \to X$ with  each $(X_n, \s)$  being a non-wandering SFT, then $(X,\s)$ is chain recurrent. \et

\bo Since $X_n \to X$, for given $\ep > 0$, there exists $N >0$ such that $X_n \subset X_\ep$ for $n > N$.

Thus for $x \in X$, there exists $x_n \in X_n$ such that for $n > N$, $x_n$ is in $\ep-$neighbourhoods of $x$ and $x_n \to x$.

Thus there exist sequence $\{k_n\} \nearrow \infty$ such that

$$   x_{[-k_n, k_n ]} = {x_n}_{[-k_n, k_n]}.$$

Since each $X_n$ is non-wandering, for  $x_n \in X_n$  there exists $M_n >0$ and $z_n \in X_n$ such that

$$ {z_n}_{[-k,k]} = {x_n}_{[-k,k]} \Longrightarrow {z_n}_{[M_n -k, M_n +k]} = {x_n}_{[-k,k]}$$

\bigskip

Choose $m > N$ such that $k_m > k$, then there exists    $L_1, \ldots, L_t > 0$ and $z_1, \ldots, z_t \in X$ such that
${z_1}_{[-k,k]} = x_{[-k,k]}, {z_1}_{[L_1-k,L_1+k]} = {z_2}_{[-k,k]}, \ \ldots, \ \text{and} \ {z_{t-1}}_{[L_{t-1}-k, L_{t-1}+k]} = {z_t}_{[-k,k]} \ \text{with } \ {z_t}_{[L_t-k, L_t+k]} = x_{[-k,k]}.$

This proves the chain recurrence of $X$.\eo

\bt For   $X \in \fS$, if $X_n \to X$ with  each $(X_n, \s)$  a transitive(irreducible) SFT, then $(X,\s)$ is chain transitive. \et

\bo Since $X_n \to X$, for given $\ep > 0$, there exists $N >0$ such that $X_n \subset X_\ep$ for $n > N$.

Thus for $x,y \in X$, there exists $x_n,y_n \in X_n$ such that for $n > N$, $x_n, y_n$ are in $\ep-$neighbourhoods of $x,y$ and $x_n \to x$ and $y_n \to y$.

Thus there exist sequence $\{k_n\} \nearrow \infty$ such that

$$   x_{[-k_n, k_n ]} = {x_n}_{[-k_n, k_n]} \ \ \& \ \ y_{[-k_n, k_n ]} = {y_n}_{[-k_n, k_n]}.$$

Since each $X_n$ is transitive, for every $x_n,y_n \in X_n$  there exists $M_n >0$ and $z_n \in X_n$ such that

$$ {z_n}_{[-k,k]} = {x_n}_{[-k,k]} \Longrightarrow {z_n}_{[M_n -k, M_n +k]} = {y_n}_{[-k,k]}$$

\bigskip

Choose $m > N$ such that $k_m > k$, then there exists    $L_1, \ldots, L_t > 0$ and $z_1, \ldots, z_t \in X$ such that
$${z_1}_{[-k,k]} = x_{[-k,k]}, {z_1}_{[L_1-k,L_1+k]} = {z_2}_{[-k,k]}, \ \ldots, \ \text{and} \ {z_{t-1}}_{[L_{t-1}-k, L_{t-1}+k]} = {z_t}_{[-k,k]} \ \text{with } \ {z_t}_{[L_t-k, L_t+k]} = y_{[-k,k]}.$$

This proves the chain transitivity of $X$.\eo

And the discussions in Remark \ref{ct=t} gives,

\bc If  $X_n \to X$ with  each $(X_n, \s)$  a transitive(irreducible) SFT and $X \in \sft$, then $(X,\s)$ is transitive.\ec

\bex \label{ttoct} We note that the limit of transitive SFTs need not be transitive.

\smallskip

We consider  SFTs $(X_n, \s)  \ \forall \ n \in \N$,  with alphabet $\{1,2,3\}$ and forbidden words

\noindent
$\cF_n = \{21^k2 , 21^k3, 22, 23, 32, 33, 31^k2, 31^k3: k \leq n \}$ respectively.

We note that in each $X_n$, a $2$ or $3$ must be followed by a block of at least $n+1$ $1s$.

Also $X_1 \subsetneqq X_2 \subsetneqq \ldots$ and $X_n \to X$ where $$X = \bigcap \limits_{n \in \N} X_n = \{1^{\infty}\} \cup \{\s^n(1^{\infty}.21^{\infty}): n \in \Z\} \cup \{\s^n(1^{\infty}.31^{\infty}): n \in \Z\}.$$

Each $(X_n, \s)$ is transitive but $(X, \s)$ is not transitive.

\medskip

We note that $(X,\s)$ is chain transitive.\eex

\bex \label{even} We note that the limit of transitive SFTs can be transitive even when not SFT.

\smallskip

Let $(X,\s)$ be the \emph{even shift} with alphabet $\{1,2\}$ and forbidden words
$\cF = \{21^k2 : k \in 2\N +1 \}$.

We consider  SFTs $(X_n, \s)  \ \forall \ n \in \N$,  with alphabet $\{1,2\}$ and forbidden words
$\cF_n = \{21^k2 : k \in 2\N +1 \text{and} \ k \leq n \}$ respectively.

We note that in each $X_n$, a $2$  must be followed by a block of at least $n+1$ $1s$.

Also $X_1 \supsetneqq X_2 \supsetneqq \ldots$ and $X_n \to X$.

\medskip

We note that each $(X_n, \s)$ is transitive and $(X, \s)$ is also transitive.
\eex

\bt For   $X \in \fS$, if $X_n \to X$ with  each $(X_n, \s)$  a mixing SFT, then $(X,\s)$ is chain mixing. \et

\bo Since $X_n \to X$, for given $\ep > 0$, there exists $N >0$ such that $X_n \subset X_\ep$ for $n > N$.

Thus for $x,y \in X$, there exists $x_n,y_n \in X_n$ such that for $n > N$, $x_n, y_n$ are in $\ep-$neighbourhoods of $x,y$ and $x_n \to x$ and $y_n \to y$.

Thus there exist sequence $\{k_n\} \nearrow \infty$ such that

$$   x_{[-k_n, k_n ]} = {x_n}_{[-k_n, k_n]} \ \ \& \ \ y_{[-k_n, k_n ]} = {y_n}_{[-k_n, k_n]}.$$

Since each $X_n$ is mixing, for every $x_n,y_n \in X_n$  there exists $M_n >0$ and $z_{n_t} \in X_n$ for all $N_t \geq M_n$ such that

$$ {z_{n_t}}_{[-k,k]} = {x_n}_{[-k,k]} \Longrightarrow {z_{n_t}}_{[N_t -k, N_t +k]} = {y_n}_{[-k,k]}$$

\bigskip

Choose $m > N$ such that $k_m > k$, then there exists $M > 0$,  and $w_1, \ldots, w_j \in X$ and $N_1, \ldots, N_j > 0$ for $j \geq M$, such that

\begin{center}
	${w_1}_{[-k,k]} = x_{[-k,k]}$
	
	${w_2}_{[-k, k]} = {w_1}_{[N_1-k,N_1 +k]}$
	
	$ \ldots $
	
	${w_j}_{[N_j-k, N_j+k]} = y_{[-k,k]}.$ \end{center}

This proves the chain mixing of $X$.\eo

And the discussions in Remark \ref{ct=t} gives,

\bc If  $X_n \to X$ with  each $(X_n, \s)$  a mixing SFT and $X \in \sft$, then $(X,\s)$ is mixing.\ec

\br We note that Example \ref{ttoct} gives an example of mixing SFTs whose limit is chain mixing but not mixing.  Example \ref{even} gives an example of mixing SFTs whose limit is mixing though not SFT. Also Example \ref{ttoct} and Example \ref{even} both illustrate Theorem \ref{dicho}.\er

We can say something more here. Let $X_n \to X$ and let $x,y \in X$. For $\ep > 0$ let $B_\ep(x), \ B_\ep(y)$ denote the $\ep$-balls centered around $x,y$ respectively in $\N^{\Z}$.

If $(X,\s)$ is also transitive then we note that the hitting times
$$N(B_\ep(x) \cap X, B_\ep(y) \cap X) = \lim \limits_{n \to \infty} N(B_\ep(x) \cap X_n, B_\ep(y) \cap X_n) $$

where this limit is taken in $2^{\N}$.

This is simple to observe since for every $\ep > 0$, and $x,y \in X$, there exists an $N > 0$ for which $B_\ep(x) \cap X_n \neq \emptyset$, and $B_\ep(y) \cap X_n \neq \emptyset$ for $n \geq N$.

\subsection{Convergence of Dynamics for Polygonal Billiards}

The results that we have established in Theorems \ref{ideal2}, \ref{compact2},   \ref{semi-ideal2} indicate the coming together of geometric and dynamical convergence  for the class of semi-ideal rational polygons in a sense that we discuss ahead. We construct a sequence of polygons $\Pi^n$ in the semi-ideal rational class that converges to a polygon $\Pi$ in the same class. Then, our results indicate that the corresponding sequence of compactifications of the space of codes $\tilde{X}^n$ must converge to the compactification of the space of codes for $\Pi$.

For a sequence of polygons $\{\Pi^n \}_{n \in \mathbb{N}}$
and a polygon $\Pi$ with $k$ vertices in $\mathbb{D}$, we  label the vertices of each $\Pi^n$ by $v^n_i$ with $i \in \{1,2,...,k\} $ and the vertices of $\Pi$ by $v_i$ with   $i \in \{1,2,...,k\} $. If $v^n_i \rightarrow v_i$ for each $i$, then we call it the \emph{convergence of a polygonal sequence} and denote it as $\Pi^n \rightarrow \Pi$. This convergence is in accordance with  the \emph{Gromov-Hausdorff metric}  as given in Equation (\ref{2}).

\bigskip

 We consider the full shift  on $k$ symbols and denote the space of all compact subsets of $\{1,2,...,k\}^{\mathbb{Z}}$ by $\M$. In view of the above discussion, we talk about the convergence of a sequence of compactifications of the spaces of codes of the corresponding polygonal billiard tables under the ambit of the metric structure provided by $(\M, d_{GH})$. It is to be noted here that for a polygonal billiard table of above mentioned class, the corresponding compactification of the space of codes uniquely determines the space of codes itself. This will allow us to trace back our sequence of polygons starting with the convergence of a sequence of associated compactifications of the space of codes in a sense that we describe ahead.

\bt \label{convergence} Let $\{\Pi^n\}_{n \in \mathbb{N}}$ be a sequence of semi-ideal rational polygons in $\mathbb{D}$ each having  $k$ vertices and $\Pi$ be another such polygon. Further, suppose $X^n$ be the corresponding space of codes for $\Pi^n$  and $\tilde{X}^n$ the corresponding compactification, for each $n \in \mathbb{N}$. Let $X$ be the space of codes for $\Pi$ with compactification $\tilde{X}$.  Then $\Pi^n \rightarrow \Pi $ implies  $\tilde{X^n} \rightarrow \tilde{X}$. 
\et

\bo
Let $\Pi^n \rightarrow \Pi$ then $v^n_i \rightarrow v_i\ \forall i\ \in \{1,2,...,k\} $. We note here that for the given counter-clockwise labelling of each $\Pi^n$, we label $\Pi$ as suggested by the above convergence. Consider any arbitrary $j \in \{1,2,...,k\}$. If $\Omega_j$ is an ideal vertex, then with $\Omega^n_j \rightarrow 0$, we get $\forall\ \epsilon >0,\ |\Omega^n_j - 0| < \epsilon $ for sufficiently large $n$. Now, moving to a subsequence with non-zero terms, if needed, we get $\lambda^n_j = \pi/ \Omega^n_j > 1/ \epsilon$ for sufficiently large $n$. Note that the zero terms of the sequence, if present satisfy the below mentioned criteria vacuously. Thus, for sufficiently large $n$, we have $(j,j+1)^{\lambda^n_j + 1} \in (\cL(X^n))^c$. Therefore, as $n \rightarrow \infty ,\ (\cL(X^n))^c$ contains $(j,j+1)^p\ \forall\ p \in \mathbb{N}$.

\bigskip

If $\Omega_j $ is a non-zero rational, then $\forall\  \epsilon > 0,\ |\Omega^n_j - \Omega_j| < \epsilon $ for sufficiently large $n$. Say, $\Omega_j = \pi / m$, then if $\epsilon < \pi/m(m+1)$, we get $\Omega^n_j = \Omega_j$ for sufficiently large $n$. Thus, $\lambda^n_j = \lambda_j$ for sufficiently large $n$, which further implies that $(j,j+1)^{\lambda_j} \in (\cL(X^n))^c$ for sufficiently large $n$. Therefore, $\tilde{X}^n \rightarrow \tilde{X}$ in $d_{GH}$ metric. 
\eo

\br We note that we can recover $X_n$ and $X$ from $\tilde{X}_n$ and $\tilde{X}$ to see that $X_n \rightarrow X$ in the sense of Gromov. \er

\br The converse of the Theorem \ref{convergence} is more subtle. We describe the case of an ideal triangle $\Pi$. We can assume that $\Pi $ is so placed that it contains the center of $\mathbb{D}$. This can be done as the coding rules are not dependent on the position of the vertices, thereby the polygon can be tweaked as required.
Its corresponding $\tilde{X}$ is described by the forbidden set $\cF_{\tilde{X}} = \{11,22,33\}$. Let us pick up a sequence $\tilde{X^n}$ given by the corresponding forbidden sets  $\cF_{\tilde{X^n}} = \{11,22,33,(12)^{n+3}\}$. Then $\tilde{X^n} \rightarrow \tilde{X}$ in the $d_{GH}$ metric. For $\tilde{X}$ the corresponding $X$ is determined by the coding rules of Theorem \ref{ideal2}. We call the center of $\mathbb{D}$ as $O$ and construct the radial Euclidean lines $Ov_1, Ov_2, Ov_3$. Next, we choose a sequence of polygons $\Pi^n$ by taking $v^n_1 \equiv v_1, v^n_2 \equiv v_2 $ and placing $v^n_3$ on $Ov_3$. Note that for each $n$ the choice of  $v^n_3$ is unique under the requirement that $\Omega^n_3 = \pi / (n+2)$. This ensures that $\Pi^n \rightarrow \Pi$ and gives us a partial answer to the converse of Theorem \ref{convergence}. The same can be said about an arbitrary ideal polygon. We also remark here that the above construction still holds good for any semi-ideal rational polygon with the restriction that exactly one vertex is in $\mathbb{D}$ with all other vertices sitting on $\partial \mathbb{D}$. \er

\section*{Acknowledgements}

We thank \textsc{Mike Boyle} for many helpful discussions on subshifts, and a very generous \textsc{Anonymous Referee} for many useful suggestions improving the readability of this article.

\end{document}